\documentclass[11pt,oneside]{article}

\usepackage[square,numbers]{natbib}

\usepackage{amsmath, amssymb}
\usepackage{calrsfs}
\DeclareMathAlphabet{\pazocal}{OMS}{zplm}{m}{n}
\usepackage[shortlabels]{enumitem}
\usepackage{xcolor}
\usepackage{amsthm}
\usepackage[nottoc]{tocbibind}
\usepackage{tocloft}
\usepackage{bbm}					
\numberwithin{equation}{section}
 
\newtheorem{thm}{Theorem}[section]
\newtheorem{lem}{Lemma}[section]
\newtheorem{prop}{Proposition}[section]
\newtheorem{cor}{Corollary}[section]

\theoremstyle{definition}
\newtheorem{defn}{Definition}[section]

\theoremstyle{remark}

\theoremstyle{remark}
\newtheorem{rem}{Remark}[section]

\theoremstyle{theorem}
\newtheorem{assump}{Assumption}[section]

\title{Concentration of Measure for the Graphon Particle System}

\vspace{2mm}

\author{  
\textsc{Erhan Bayraktar} 
\thanks{
E. Bayraktar is supported in part by the National Science Foundation under grant DMS-2106556 and by the Susan M. Smith Professorship.
}  
\and
\textsc{Donghan Kim} 
}

\begin{document}

\maketitle


\begin{abstract}
\noindent
We study heterogeneously interacting diffusive particle systems with mean-field type interaction characterized by an underlying graphon and their finite particle approximations. Under suitable conditions, we obtain exponential concentration estimates over a finite time horizon for both 1 and 2 Wasserstein distances between the empirical measures of the finite particle systems and the averaged law of the graphon system, extending the work of Bayraktar-Wu~\cite{BW}.
\end{abstract}


\let\thefootnote\relax\footnotetext{
	{\it MSC 2020 subject classifications:} 05C80, 60J60, 60K35.
}
\let\thefootnote\relax\footnotetext{
	{\it Keywords and phrases:} graphon, graphon particle system, heterogeneous interaction, mean-field interaction, networks, concentration bounds, transport inequalities, $L^1$-Fourier class
}

\smallskip

\tableofcontents

\input amssym.def
\input amssym

\section{Introduction}
\label{sec: introduction}

In this article, we study concentration of measures related to the graphon particle system and its finite particle approximations. This work is a continuation of earlier papers \cite{BCW, BW, BW2}. A graphon particle system consists of uncountably many heterogeneous particles $X_u$ for $u \in [0, 1]$ with which their interactions are characterized by a graphon. More precisely, for a fixed $T > 0$ and $d \in \mathbb{N}$, we consider the following system 
\begin{align}
X_u(t) & = X_u(0) 		\nonumber
+ \int_0^t \bigg[ \int_0^1 \int_{\mathbb{R}^d} \phi \big(X_u(s), y \big) G(u, v) \, \mu_{v, s} (dy) \, dv + \psi \big(X_u(s)\big) \bigg] ds 
+ \sigma B_u(t),
\\ 
& \mu_{u, t} \text{ is the probability distribution of } X_u(t), \text{ for every } u \in [0, 1] \text{ and } t \in [0, T],		\label{eq : Graphon}
\end{align}
where $\{B_u\}_{u \in [0, 1]}$ is a family of i.i.d. $d$-dimensional Brownian motions, $\{X_u(0)\}_{u \in [0, 1]}$ is a collection of independent (but not necessarily identically distributed) $\mathbb{R}^d$-valued random variables with law $\mu_u(0)$, independent of $\{B_u\}_{u \in [0, 1]}$ for each $u \in [0, 1]$, defined on a filtered probability space $(\Omega, \mathcal{F}, \{\mathcal{F}_t\}, \mathbb{P})$. Two functions $\phi : \mathbb{R}^d \times \mathbb{R}^d \rightarrow \mathbb{R}^d$ and $\psi : \mathbb{R}^d \rightarrow \mathbb{R}^d$ represent pairwise interactions between the particles and single particle drift, respectively. $\sigma \in \mathbb{R}^{d \times d}$ is a constant and $G : [0, 1] \times [0, 1] \rightarrow [0, 1]$ is a graphon, that is, a symmetric measurable function.

\smallskip

Along with the graphon particle system, we introduce two finite particle systems with heterogeneous interactions, which approximate \eqref{eq : Graphon}. For a fixed, arbitrary $n \in \mathbb{N}$ and each $i \in [n] := \{ 1, \cdots, n\}$, we first consider the not-so-dense analogue of \eqref{eq : Graphon} introduced in Section~4 of \cite{BCW}:
\begin{equation}	\label{eq : SDE}
X^n_i(t) = X_{\frac{i}{n}}(0) + \int_0^t \Big[ \frac{1}{n p(n)} \sum_{j=1}^n \xi^n_{ij} \phi \big(X^n_i(s), X^n_j(s) \big) + \psi \big( X^n_i(s) \big) \Big] ds + \sigma B_{\frac{i}{n}}(t),
\end{equation}
where $\{p(n)\}_{n \in \mathbb{N}} \subset (0, 1]$ is a sequence of numbers and $\{\xi^n_{ij}\}_{1 \le i, j \le n}$ are independent Bernoulli random variables satisfying
\begin{equation*}
\xi^n_{ij} = \xi^n_{ji}, \qquad \mathbb{P}(\xi^n_{ij} = 1) = p(n) G(\frac{i}{n}, \, \frac{j}{n}), \qquad \text{for every } i, j \in [n],
\end{equation*}
independent of $\{B_{i/n}, X_{i/n}(0) : i \in [n]\}$. Here, $p(n)$ represents the global sparsity parameter and the strength of interaction between the particles in \eqref{eq : SDE} is scaled by $np(n)$, the order of the number of neighbors, as in mean-field systems on Erd\"os-R\'enyi random graphs \cite{BHAMIDI20192174, MFG_Delarue, Oliveria:Reis}; the convergence of $p(n) \rightarrow 0$ as $n \rightarrow \infty$ implies that the graph is sparse, but we shall consider the not-so-dense case $np(n) \rightarrow \infty$ of diverging average degree in the random graph.

\smallskip

The other finite particle approximation system is given by
\begin{equation}	\label{eq : SDE bar}
\bar{X}^n_i(t) = X_{\frac{i}{n}}(0) + \int_0^t \Big[ \frac{1}{n} \sum_{j=1}^n \phi \big(\bar{X}^n_i(s), \bar{X}^n_j(s) \big) G(\frac{i}{n}, \frac{j}{n}) + \psi \big( \bar{X}^n_i(s) \big) \Big] ds + \sigma B_{\frac{i}{n}}(t).
\end{equation}
Since this system has a nonrandom coefficient for the interaction term (but still models heterogeneous interaction via the graphon), it is easier to analyze than the other finite particle system \eqref{eq : SDE}. We note that three systems \eqref{eq : Graphon} - \eqref{eq : SDE bar} are coupled in a sense that they share initial particle locations $X_{i/n}(0)$ and Brownian motions $B_{i/n}$ for $i \in [n]$.

\smallskip

Law of large numbers~(LLN)-type of convergence results for the systems \eqref{eq : SDE} and \eqref{eq : SDE bar} to the graphon particle system \eqref{eq : Graphon} under suitable conditions are studied in \cite{BCW}. The exponential ergodicity of the two systems \eqref{eq : Graphon} and \eqref{eq : SDE}, as well as the uniform-in-time convergence of \eqref{eq : SDE} to \eqref{eq : Graphon} under a certain dissipativity condition, are presented in \cite{BW2}. There are more studies of the graphon particle systems \cite{BW, Coppini:Graphon}, and works of associated heterogeneously interacting finite particle models \cite{BCN, Coppini, Delattre2016ANO, Lacker2021ACS, Lucon, Oliveria:Reis}. These studies are recently arisen since graphons have been widely applied in mean-field game theory for both static and dynamic cases, see e.g. \cite{Carmona:SGG2, BWZ, Caines:Huang, Caines:Huang2, Carmona:Cooney, GCH, Gao:Tchuendom:Caines, Parise2019GraphonGA, TCH2, TCH, VMV} and references therein.

\smallskip

Among these studies, our work is particularly linked to \cite{BW}. Denoting $W_1$ the 1-Wasserstein distance and defining the empirical measures of the three particle systems at time $t \in [0, T]$
\begin{equation}	\label{def : empirical measures at t}
L_{n, t} := \frac{1}{n} \sum_{i=1}^n \delta_{(X^n_i(t))}, \qquad 
\bar{L}_{n, t} := \frac{1}{n} \sum_{i=1}^n \delta_{(\bar{X}^n_i(t))}, \qquad
\widetilde{L}_{n, t} := \frac{1}{n} \sum_{i=1}^n \delta_{(X_{i/n}(t))},
\end{equation}
along with the averaged law $\widetilde{\mu}_t := \int_0^1 \mu_{u, t} \, du$ of the graphon system \eqref{eq : Graphon}, concentration bounds of the types $\mathbb{P} \big[ \sup_{0 \le t \le T} W_1(\bar{L}_{n, t}, \widetilde{\mu}_t) > \epsilon \big]$, $\sup_{t \ge 0} \mathbb{P} \big[  W_1(\bar{L}_{n, t}, \widetilde{\mu}_t) > \epsilon \big]$ for $\epsilon > 0$ are computed in \cite{BW} under certain conditions. In particular, uniform-in-time concentration bound of the latter type is studied in an infinite time horizon setting under an extra dissipativity condition on $\psi$. These results are established by computing certain sub-Gaussian estimates rather directly with the moment generating function of the standard normal random vector~(Lemmas~3.7 - 3.10 of \cite{BW}). 

\smallskip

In contrast, our work focuses on the case of finite time horizon and deals with more general sparsity sequence $\{p(n)\}_{n \in \mathbb{N}} \subset (0, 1]$ for \eqref{eq : SDE}, whereas the results of \cite{BW} only cover the dense graphs, i.e., $p(n) \equiv 1$. Our argument adopts the method of \cite{DLR:AOP} as follows. When $\bar{X}^n = (\bar{X}^n_1, \cdots, \bar{X}^n_n)$ represents the state of the so-called Nash equilibrium of a symmetric $n$-player stochastic differential game and $\widetilde{\mu}$ is the measure flow of the unique equilibrium of the corresponding mean-field game, the authors of \cite{DLR:AOP} compute concentration bounds of the probabilities
\begin{equation}	\label{concentration bounds 1, 2}
	\mathbb{P} \big[ \sup_{0 \le t \le T} W_p(\bar{L}_{n, t}, \widetilde{\mu}_t) > \epsilon \big], \qquad p \in \{1, 2\},
\end{equation}
with the notation of \eqref{def : empirical measures at t}. Their argument uses transportation inequalities in \cite{Djellout:Guillin:Wu} to show that Lipschitz functions of $\bar{X}^n$ concentrate around their means, and obtains the aforementioned bounds from this concentration property. We apply a similar approach to the finite particle system \eqref{eq : SDE bar} to obtain the bound of the probability that Lipschitz function values of the particles on the space $\big(C([0, T] : \mathbb{R}^d)\big)^n$, deviate from their means in Theorems~\ref{thm : concen ineq X bar1} and \ref{thm : concen ineq X bar2}. Combining this bound with the facts in Section~\ref{subsec : expectation} that the expectations of the $W_2$-distances between the empirical measures in \eqref{def : empirical measures at t} converge to zero as the number of particles goes to infinity, we show the concentration result of \eqref{concentration bounds 1, 2}. In particular, obtaining the same exponential bound in $n$ for the probability \eqref{concentration bounds 1, 2} in terms of $W_2$-metric in Theorem~\ref{thm : W2 L bar}, as in the case of $p = 1$ in Theorem~\ref{thm : W1 L bar}, is the new result.

\smallskip

Moreover, inspired by the proof in \cite{Oliveria:Reis}, we compare the particles $X^n_i$ and $\bar{X}^n_i$ to improve the exponential bound in \cite{BW} for $\mathbb{P} \big[ \sup_{0 \le t \le T} W_p(L_{n, t}, \widetilde{\mu}_t) > \epsilon \big]$ when $p = 1$, at the expense of an assumption on the interaction function $\phi$, namely being a member of the $L^1$-Fourier class. When $p=2$, we also present the similar exponential bound for the system \eqref{eq : SDE} on the dense graphs~($p(n) \equiv 1$). Without such condition on $\phi$, we have the same bound, but in terms of the bounded Lipschitz metric~($d_{BL}$-metric), a weaker metric than $W_1$. 

\smallskip

This paper is organized as follows. In Section~\ref{sec: setting} we introduce the notation, state the assumptions, and recall some of the relevant existing results concerning the particle systems \eqref{eq : Graphon} - \eqref{eq : SDE bar} and the other preliminary results. Section~\ref{sec: concentration} provides our main results and Section~\ref{sec : proofs} gives proofs of these results.

\bigskip

\section{Preliminaries}
\label{sec: setting}

In this section, we first introduce the notation which will be used throughout this paper. We then state several assumptions with some of the basic results on the particle systems \eqref{eq : Graphon} - \eqref{eq : SDE bar} from \cite{BCW, BW}, and provide several well-known results regarding transportation cost inequalities without proof. Finally, Bernstein's inequality with the concept of $L^1$-Fourier class will be introduced.

\medskip

\subsection{Notation}

Given a metric space $(S, d)$ and a function $f : S \rightarrow \mathbb{R}$, we define
\begin{align*}
	\vert\vert f \vert\vert_{\infty} &:= \sup \{ \vert f(x) \vert : x \in S\},
	\\
	\vert\vert f \vert\vert_{Lip} &:= \sup \Big\{ \frac{\vert f(x) - f(y) \vert}{d(x, y)} : x, y \in S, \, x \neq y \Big\},
	\\
	\vert\vert f \vert\vert_{BL} &:= 2 \big( \vert \vert f(x) \vert \vert_{\infty} + \vert\vert f \vert\vert_{Lip} \big),
\end{align*}
and we say that $f$ is Lipschitz~(bounded Lipschitz) if $\vert\vert f \vert\vert_{Lip} < \infty$~($\vert\vert f \vert\vert_{BL} < \infty$), respectively. In particular, $f$ is called $a$-Lipschitz if  $\Vert f \Vert_{Lip} = a$, and note that $\vert\vert f \vert\vert_{BL} \le 1$ implies $\vert f(x) - f(y) \vert \le \vert x-y \vert \wedge 1$.

\smallskip

Denote by $\pazocal{P}(S)$ the space of Borel probability measures on $S$ and we shall use the standard notation $\langle \mu, \varphi \rangle := \int_S \varphi \, d\mu$ for integrable functions $\varphi$ and measures $\mu$ on $S$. When $(S, \vert\vert \cdot \vert\vert)$ is a normed space, we write $\pazocal{P}^p(S, \vert\vert\cdot\vert\vert)$ for the set of $\mu \in \pazocal{P}(S)$ satisfying $\langle \mu, \vert\vert\cdot\vert\vert^p \rangle < \infty$ for a given $p \in [1, \infty)$. We denote by $Lip(S, \vert\vert\cdot\vert\vert)$ the set of $1$-Lipschitz functions, i.e., $f : S \rightarrow \mathbb{R}$ satisfying $\vert f(x) - f(y) \vert \le \vert\vert x-y \vert\vert$ for every $x, y \in S$. 

\smallskip

For a separable Banach space $(S, \vert\vert\cdot\vert\vert)$, we endow $\pazocal{P}^p(S, \vert\vert\cdot\vert\vert)$ with the $p$-Wasserstein metric
\begin{equation*}
	W_{p, (S, \vert\vert\cdot\vert\vert)} (\mu, \nu) := \inf_{\pi} \Big(\int_{S \times S} \vert\vert x-y \vert\vert^p \pi(dx, \, dy) \Big)^{1/p},  \qquad p \ge 1, 
\end{equation*}
where the infimum is taken over all probability measures $\pi$ on $S \times S$ with first and second marginals $\mu$ and $\nu$. We also write the product space $S^n := S \times \cdots \times S$, equipped with the $\ell^p$ norm for any $p \ge 1$
\begin{equation*}
	\vert\vert x \vert\vert_{n, p} = \Big( \sum_{i=1}^n \vert\vert x_i \vert\vert^p \Big)^{1/p},
\end{equation*}
for $x = (x_1, \cdots, x_n) \in S^n$. When the space $S$ or the norm $\vert\vert\cdot\vert\vert$ is understood, we sometimes omit it from the above notations. 

\smallskip

Denote by $C([0, T] : S)$ the space of continuous functions from $[0, T]$ to $S$, and $\vert\vert x \vert \vert_{\star, t} := \sup_{0 \le s \le t}\vert x_s \vert$, where $\vert \cdot \vert$ is the usual Euclidean norm on $\mathbb{R}^d$ for $x \in C([0, T] : \mathbb{R}^d)$ and $t \in [0, T]$. We write $\pazocal{L}(X)$ the probability law of a random variable $X$ and $[n] := \{1, \cdots, n\}$ for any $n \in \mathbb{N}$. We use $K$ to denote various positive constants throughout the paper and its value may change from line to line. 

\smallskip

For a Polish space $(S, d)$ with Borel $\sigma$-field $\pazocal{S}$, we also consider the space of probability measures over $(S, \pazocal{S})$ endowed with the topology of weak convergence, which is metrized by the $BL$-metric, defined for $\mu, \nu \in \pazocal{P}(S)$ as
\begin{equation}		\label{def : BL metric}
	d_{BL}(\mu, \nu) := \sup \bigg\{ \Big\vert \int_S \, f \, d(\mu-\nu) \Big\vert : f:S \rightarrow \mathbb{R} \text{ with } \vert\vert f \vert\vert_{BL} \le 1 \bigg\}.
\end{equation}
Note the dual representation of the $1$-Wasserstein metric
\begin{equation}		\label{def : W1 metric}
	W_1(\mu, \nu) :=  \sup \bigg\{ \Big\vert \int_S \, f \, d(\mu-\nu) \Big\vert : f:S \rightarrow \mathbb{R} \text{ with } \vert\vert f \vert\vert_{Lip} \le 1 \bigg\},
\end{equation}
along with the relationship $d_{BL} \le W_1$. We shall also use the notation for given $\mu, \nu \in \pazocal{P}(C([0, T] : \mathbb{R}^d))$
\begin{equation*}
	W_{p, t}(\mu, \nu) := \inf_{\pi} \Big( \int \vert \vert x- y \vert \vert^p_{\star, t} \, \pi(dx, \, dy) \Big)^{1/p}, \qquad t \in [0, T], \quad p \ge 1,
\end{equation*}
where the infimum is taken over all probability measures $\pi$ with marginals $\mu$ and $\nu$.

\smallskip

Let us define three $n \times n$ random matrices $P^{(n)}$, $\bar{P}^{(n)}$, and $D^{(n)}$, concerning the systems \eqref{eq : SDE}, \eqref{eq : SDE bar} for every $n \in \mathbb{N}$ with entries
\begin{align}
P^{(n)}_{i, j} &:= \frac{\xi^n_{ij}}{np(n)}, \qquad \qquad i, j \in [n],	\nonumber
\\
\bar{P}^{(n)}_{i, j} &:= \frac{1}{n}G(\frac{i}{n}, \frac{j}{n}), ~\qquad i, j \in [n],			\label{def : P bar}
\\
D^{(n)} &:= P^{(n)}-\bar{P}^{(n)}.		\nonumber
\end{align}
For these matrices, we define the $\ell_{\infty} \rightarrow \ell_1$ norm of an $n \times n$ matrix $A$
\begin{equation}
\vert\vert A \vert\vert_{\infty \rightarrow 1} := \sup \Big\{ \big\langle \mathbf{x}, \, A \mathbf{y} \big\rangle \, : \, \mathbf{x}, \, \mathbf{y} \in [-1, 1]^n \Big\}.
\end{equation}
This norm is known to be equivalent to the so-called \textit{cut norm}~(see (3.3) of \cite{GV}).

\smallskip

We denote the empirical measures of the approximation systems for each $n \in \mathbb{N}$
\begin{equation}	\label{def : empirical measures}
L_n := \frac{1}{n} \sum_{i=1}^n \delta_{(X^n_i)}, \qquad
\bar{L}_n := \frac{1}{n} \sum_{i=1}^n \delta_{(\bar{X}^n_i)},
\end{equation}
all of which are random elements of $\pazocal{P}(C([0, T] : \mathbb{R}^d))$.

\smallskip

We conclude this subsection by recalling the relative entropy of two probability measures $\mu, \nu$ over the same measurable space 
\begin{equation}	\label{def : entropy}
H(\mu \vert \nu) := 
\begin{cases}
\int \log \big(\frac{d\mu}{d\nu} \big) d\mu, \qquad \text{if } \mu \ll \nu;
\\ \qquad \infty, \qquad \qquad ~ \text{otherwise}.
\end{cases}
\end{equation}

\medskip

\subsection{Existence and uniqueness of the solutions}

We state the existence and uniqueness of strong solutions to the systems \eqref{eq : Graphon} - \eqref{eq : SDE bar}.

\begin{assump}	\label{assump : Lipschitz}
	\hspace{1cm}
	\begin{enumerate} [(a)]
		\item $\phi$ is bounded; $\phi$ and $\psi$ are Lipschitz, i.e.,
		there exists a constant $K > 0$ such that
		\begin{equation*}
			\big\vert \phi(x_1, y_1) - \phi(x_2, y_2) \big\vert + \big\vert \psi(x_1) - \psi(x_2) \big\vert \le K \big( \vert x_1-x_2 \vert + \vert y_1 - y_2 \vert \big)
		\end{equation*}
		holds. Moreover, the initial particles have finite second moments, i.e., 
		\begin{equation}	\label{con : finite second initial}
		\sup_{u \in [0, 1]} \mathbb{E} \big\vert X_u(0) \big\vert^2 < \infty.
		\end{equation}
		\item The map $[0, 1] \ni u \mapsto \mu_u(0) = \mathcal{L}(X_u(0)) \in \pazocal{P}(\mathbb{R}^d)$ is measurable.
	\end{enumerate}
\end{assump}

\smallskip

\begin{lem}	[The existence and uniqueness of the particle systems]\label{lem : existence of solutions}
	\hspace{1cm}
	\begin{enumerate} [(a)]
		\item Under Assumption~\ref{assump : Lipschitz} (a), two systems \eqref{eq : SDE}, \eqref{eq : SDE bar} have unique strong solutions.
		\item Under Assumption~\ref{assump : Lipschitz} (a) and (b), the graphon system \eqref{eq : Graphon} has a unique strong solution, and the map $[0, 1] \ni u \mapsto \mu_u \in \pazocal{P}\big(C([0, T] : \mathbb{R}^d)\big)$ is measurable.
	\end{enumerate}	
\end{lem}

\smallskip

Proof of Lemma~\ref{lem : existence of solutions} (a) is classical (see e.g. Theorem 5.2.9 of \cite{KS1}). Part (b) follows from Proposition 2.1 of \cite{BCW}. As pointed out in Remark~2.2 of \cite{BCW}, we note that the boundedness condition on $\phi$ in Assumption~\ref{assump : Lipschitz}(a) can be removed throughout this paper at the expense of a stronger condition $\sup_{u \in [0, 1]} \mathbb{E} \vert X_u(0) \vert^{2 + \epsilon} < \infty$ for some $\epsilon > 0$ than \eqref{con : finite second initial}. We occasionally need an even stronger condition on the initial particles as in the following.

\begin{assump}	\label{assump : independent initial}
	The initial particles $\{X_{u}(0)\}_{u \in [0, 1]}$ are independent with law $\mu_{u, 0} \in \pazocal{P}(\mathbb{R}^d)$, which satisfies
	\begin{equation}	\label{con : indep initial with exponential}
		\sup_{u \in [0, 1]} \int_{\mathbb{R}^d} e^{\kappa \vert x \vert^2} \mu_{u, 0} (dx) < \infty, \qquad \text{for some } \kappa > 0.
	\end{equation}
\end{assump}

\smallskip

We observe later that the condition \eqref{con : indep initial with exponential} is equivalent to \eqref{con : initial wasser1} from Lemma~\ref{lem : Wasser equivalence}. Under this stronger assumption, we have the finite fourth moment of the solution to \eqref{eq : Graphon}. The proof is standard, hence is omitted~(see, e.g. \cite{Sznitman}, or Proposition~2.1 of \cite{BW}).

\smallskip

\begin{lem}	\label{lem : fourth moment}
	Under Assumptions~\ref{assump : Lipschitz}, \ref{assump : independent initial}, the solution to \eqref{eq : Graphon} satisfies
	\begin{equation*}
		\sup_{u \in [0, 1]} \sup_{t \in [0, T]} \mathbb{E} \big[ \big\vert X_{u}(t) \big\vert^4 \big] < \infty.
	\end{equation*}
\end{lem}

\medskip

\subsection{Continuity of the graphon system}

The following result, which states the continuity of the graphon system \eqref{eq : Graphon}, is from Theorem~2.1 of \cite{BCW}.

\begin{assump}	\label{assump : continuity}
	There exists a finite collection of subintervals $\{I_i : i \in [N]\}$ for some $N \in \mathbb{N}$, satisfying $\cup_{i=1}^N I_i = [0, 1]$. For each $i, j \in [N]$:
	\begin{enumerate} [(a)]
		\item The map $I_i \ni u \mapsto \mu_u(0) \in \pazocal{P}(\mathbb{R}^d)$ is continuous with respect to the $W_2$-metric.
		\item For each $u \in I_i$, there exists a Lebesgue-null set $N_u \subset [0, 1]$ such that $G(u, v)$ is continuous at $(u, v) \in [0, 1] \times [0, 1]$ for each $v \in [0, 1] \setminus N_u$.
		\item There exists $K > 0$ such that
		\begin{align*}
			W_2(\mu_{u_1}(0), \mu_{u_2}(0)) &\le K \vert u_1 - u_2 \vert, \qquad \qquad \qquad \qquad u_1, u_2 \in [0, 1],
			\\
			\big\vert G(u_1, v_1) - G(u_2, v_2) \big\vert &\le K \big( \vert u_1 - u_2 \vert + \vert v_1 - v_2 \vert \big), \qquad (u_1, v_1), \, (u_2, v_2) \in I_i \times I_j.
		\end{align*}
	\end{enumerate}
\end{assump}

\smallskip

\begin{lem}	 \label{lem : continuity of graphon system}
	Suppose that Assumption~\ref{assump : Lipschitz} holds.
	\begin{enumerate} [(a)]
		\item (Continuity) Under Assumption~\ref{assump : continuity} (a) and (b), the map $I_i \ni u \mapsto \mu_u \in \pazocal{P}\big(C([0, T] : \mathbb{R}^d)\big)$ is continuous with respect to the $W_{2, T}$ metric for every $i \in [N]$.
		\item (Lipschitz continuity) Under Assumption~\ref{assump : continuity} (c), there exists $\kappa > 0$, which depends on $T$, such that $W_{2, T}(\mu_u, \mu_v) \le \kappa \vert u-v \vert$ whenever $u, v \in I_i$ for some $i \in [N]$.
	\end{enumerate}	
\end{lem}

\smallskip

In Lemma~\ref{lem : continuity of graphon system}(b), note that we have, in particular,
\begin{equation*}
	\sup_{\vert\vert f \vert\vert_{Lip} \le 1} \bigg\vert \int_{\mathbb{R}^d} f(x) \mu_{u, t} (dx) - \int_{\mathbb{R}^d} f(x) \mu_{v, t}(dx) \bigg\vert \le W_{2, t}(\mu_u, \mu_v) \le \kappa \vert u-v \vert, \qquad \forall \, t \in [0, T].
\end{equation*}

\medskip

\subsection{A law of large number of the mean-field particle system}

Besides the assumptions introduced in this section, we will need the following assumption on the sparsity parameter for the system \eqref{eq : SDE}, as briefly mentioned in Section~\ref{sec: introduction}.

\begin{assump}	\label{assump : np(n)}
	The sequence $\{p(n)\}_{n \in \mathbb{N}}$ in \eqref{eq : SDE} satisfies $np(n) \rightarrow \infty$ as $n \rightarrow \infty$.
\end{assump}

\smallskip

We introduce the following law of large number result for the mean-field particle system \eqref{eq : SDE}, which is Theorem~4.1 of \cite{BCW}. We write $\mu_u$ the law of $X_u$ in the graphon particle system \eqref{eq : Graphon} for each $u \in [0, 1]$, and define 
\begin{equation}	\label{def : widetilde mu}
	\widetilde{\mu} := \int_0^1 \mu_{u} \, du.
\end{equation}

\smallskip

\begin{lem}	\label{lem : LLN}
	Under Assumptions~\ref{assump : Lipschitz}, \ref{assump : continuity}, and \ref{assump : np(n)},
	\begin{equation*}
		L_n \rightarrow \widetilde{\mu}	\quad \text{in } \pazocal{P}\big(C([0, T] : \mathbb{R}^d)\big) \text{ in probability, as } n \rightarrow \infty.
	\end{equation*}
	Moreover, we have
	\begin{equation*}
		\frac{1}{n} \sum_{i=1}^n \mathbb{E}\vert\vert X^n_i - X_{\frac{i}{n}} \vert\vert^2_{\star, T} \rightarrow 0, \qquad \text{as } n \rightarrow \infty.
	\end{equation*}
\end{lem}

\medskip

\subsection{Transportation inequalities}	\label{subsec : Transportation}

In this subsection, we present some preliminary results regarding transportation inequalities. The first result is from Theorem~9.1 of \cite{ASU}, illustrating the transportation inequality with the uniform norm for the laws of diffusion processes.

\smallskip

\begin{lem}		\label{lem : wasser entropy}
	For a fixed $T>0$ and $k \in \mathbb{N}$, suppose that $X^x = \{X^x_t\}_{t \in [0, T]}$ is the unique strong solution of the SDE
	\begin{equation}	\label{eq : SDE form}
	dX^x_t = b(t, X^x)dt + \Sigma \, dW_t, \quad \forall \, t \in [0, T], \qquad X_0 = x \in \mathbb{R}^k,
	\end{equation}
	on a probability space $C([0, T] : \mathbb{R}^k)$ supporting a $k$-dimensional Brownian motion $W$. Here, $b : [0, T] \times C([0, T] : \mathbb{R}^k) \rightarrow \mathbb{R}^k$ satisfies for any $\xi$, $\eta \in C([0, T] : \mathbb{R}^k)$
	\begin{equation}	\label{con : b lipschitz}
	\big\vert b(t, \xi) - b(t, \eta) \big\vert \le L \sup_{0 \le s \le t} \big\vert \xi(s) - \eta(s) \big\vert = L \Vert \xi - \eta \Vert_{\star, t},  \qquad \forall \, t \in [0, T],
	\end{equation}
	for some constants $L > 0$ and $\Sigma \in \mathbb{R}^{k \times k}$. Let $P^x \in \pazocal{P}(C([0, T] : \mathbb{R}^k))$ be the law of $X^x$ for any $x \in \mathbb{R}^k$. Then for any $Q \in \pazocal{P}(C([0, T] : \mathbb{R}^k))$, we have
	\begin{equation*}
	W^2_{1, (C([0, T] : \mathbb{R}^k), \, \vert\vert\cdot\vert\vert_{k, 2})}(P^x, Q)
	\le 
	W^2_{2, (C([0, T] : \mathbb{R}^k), \, \vert\vert\cdot\vert\vert_{k, 2})}(P^x, Q) 
	\le 6e^{15L^2} H(Q \vert P^x),
	\end{equation*}
	where $H(Q\vert P)$ is the relative entropy of $Q$ with respect to $P$, defined in \eqref{def : entropy}.
\end{lem}

\smallskip

The following result~(Theorem~5.1 of \cite{DLR:AOP}) characterizes concentration of a probability measure with a transportation cost inequality and Gaussian integrability property. The equivalence between \eqref{eq : W entropy} and \eqref{eq : exponential ineq} is originally from Theorem~3.1 of \cite{Bobkov1991}, and the equivalence between \eqref{eq : W entropy} and \eqref{eq : sub-Gaussian} is due to Theorem~2.3 of \cite{Djellout:Guillin:Wu}.

\smallskip

\begin{lem}	\label{lem : Wasser equivalence}
	For a probability measure $\mu \in \pazocal{P}^1(S)$ on a separable Banach space $(S, \vert\vert\cdot\vert\vert)$, the following statements are equivalent up to a universal change in the positive constant $c$.
	\begin{enumerate} [(i)]
		\item The transportation cost inequality
		\begin{equation}	\label{eq : W entropy}
		W_{1, S}(\mu, \nu) \le \sqrt{2c H(\nu \vert \mu)}
		\end{equation}
		holds for every $\nu \in \pazocal{P}(S)$.
		
		\item For every $1$-Lipschitz function $f$ on $S$ and $\lambda \in \mathbb{R}$
		\begin{equation}	\label{eq : exponential ineq}
		\int_S e^{\lambda(f - \langle \mu, f \rangle )}d\mu \le \exp\Big(\frac{c\lambda^2}{2}\Big)
		\end{equation}
		holds.
		
		\item For every $1$-Lipschitz function $f$ on $S$ and $a > 0$
		\begin{equation}	\label{eq : mu exp}
		\mu \big(f - \langle \mu, f \rangle > a\big) \le \exp\Big(-\frac{a^2}{2c}\Big).
		\end{equation}
		
		\item $\mu$ is sub-Gaussian, i.e.,
		\begin{equation}	\label{eq : sub-Gaussian}
		\int_S e^{c\vert\vert x \vert\vert^2} \mu(dx) < \infty.
		\end{equation}
	\end{enumerate}
\end{lem}

\smallskip

The next result is well-known tensorization of transportation cost inequalities from Corollary~5 of \cite{Gozlan:Leonard}. The major difference between (i) and (ii) is that the inequality \eqref{ineq : tensor W2} is \textit{dimension-free}, i.e., the right-hand side does not depend on $n$.  

\smallskip

\begin{lem}	\label{lem : Wasser tensor}
	For each $n \in \mathbb{N}$, consider a set of probability measures $\{\mu_i\}_{i \in [n]} \subset \pazocal{P}(S)$ on a separable Banach space $(S, \Vert\cdot\Vert)$.
	\begin{itemize}
		\item [(i)] If the inequality $W_{1, S}(\mu_i, \nu) \le \sqrt{2c H(\nu \vert \mu_i)}$ holds for every $i \in [n]$ and $\nu \in \pazocal{P}^1(S)$, then 
		\begin{equation*}
		W_{1, (S^n, \, \vert\vert\cdot\vert\vert_{n, 1})}(\mu_1 \otimes \cdots \otimes \mu_n, \, \rho) 
		\le \sqrt{2nc H(\rho \vert \mu_1 \otimes \cdots \otimes \mu_n)},
		\end{equation*}
		holds for every $\rho \in \pazocal{P}^1(S^n)$.
		
		\item [(ii)] If the inequality $W_{2, S}(\mu_i, \nu) \le \sqrt{2c H(\nu \vert \mu_i)}$ holds for every $i \in [n]$ and $\nu \in \pazocal{P}^2(S)$, then 
		\begin{equation}	\label{ineq : tensor W2}
		W_{2, (S^n, \, \vert\vert\cdot\vert\vert_{n, 2})}(\mu_1 \otimes \cdots \otimes \mu_n, \, \rho)
		\le \sqrt{2c H(\rho \vert \mu_1 \otimes \cdots \otimes \mu_n)},
		\end{equation}
		holds for every $\rho \in \pazocal{P}^2(S^n)$.
	\end{itemize}
\end{lem}

\smallskip

We finally mention the following result on the Wasserstein distance of the empirical measures of independent but not necessarily identically distributed random variables. This is Lemma~A.1 of \cite{BW}, a generalization of Theorem~1 of \cite{Fournier:Guillin} where i.i.d. random variables are considered. This result will be used in proving Proposition~\ref{prop : expectation W2 L bar mu tilde}.

\smallskip

\begin{lem}	\label{lem : independent empirical measures}
	Let $\{Y_i\}_{i \in \mathbb{N}}$ be independent $\mathbb{R}^d$-valued random variables and define
	\begin{equation*}
	\nu_n := \frac{1}{n} \sum_{i=1}^n \delta_{Y_i}, \qquad
	\bar{\nu}_n := \frac{1}{n} \sum_{i=1}^n \pazocal{L}(Y_i).
	\end{equation*}
	For a fixed $p > 0$, assume that $\sup_{i \in \mathbb{N}} \mathbb{E} \vert Y_i \vert^q < \infty$ holds for some $q > p$. Then there exists a constant $K > 0$ depending only on $p$, $q$, and $d$ such that for every $n \ge 1$
	\begin{equation*}
	\mathbb{E} \big[ W^p_p(\nu_n, \bar{\nu}_n) \big]
	\le K \alpha_{p, q}(n) \bigg( \int_{\mathbb{R}^d} \vert x \vert^q \, \bar{\nu}_n (dx) \bigg)^{p/q},
	\end{equation*}
	where
	\begin{equation*}
	\alpha_{p, q}(n) :=
	\begin{cases}
	n^{-1/2} + n^{-(q-p)/q}, \qquad &\text{if } p > d/2 \text { and } q \neq 2p,
	\\
	n^{-1/2} \log(1+n) + n^{-(q-p)/q}, \qquad &\text{if } p = d/2 \text { and } q \neq 2p,
	\\
	n^{-p/d} + n^{-(q-p)/q}, \qquad &\text{if } p < d/2 \text { and } q \neq d/(d-p).
	\end{cases}
	\end{equation*}
\end{lem}

\medskip

\subsection{Bernstein's inequality and $L^1$-Fourier class}	\label{subsec: EE finite}

When comparing two approximation systems \eqref{eq : SDE} and \eqref{eq : SDE bar}, controlling the matrix $D^{(n)}$ of \eqref{def : P bar} is essential. Thus, we introduce the following concentration of $D^{(n)}$ in terms of $\Vert \cdot \Vert_{\infty \rightarrow 1}$ norm, which is from Lemma~2 of \cite{Oliveria:Reis}. Its proof is straightforward application of Bernstein's inequality~(Lemma~\ref{lem : Bernstein}, or Bennett's inequality) with the distribution of independent $n^2$ entries of the matrix $D^{(n)}$. We will use Bernstein's inequality again in Section~\ref{subsec : concentration toward GS} to prove Lemma~\ref{lem : DTD}, an elaboration of Lemma~\ref{lem : cut norm}.

\smallskip

\begin{lem}	 \label{lem : cut norm}
	For any $0 < \eta \le n$, we have
	\begin{equation*}
	\mathbb{P} \bigg[ \frac{\vert\vert D^{(n)} \vert\vert_{\infty \rightarrow 1}}{n} > \eta \bigg] 
	\le \exp \bigg( - \frac{\eta^2 n^2p(n)}{2+\frac{\eta}{3}} \bigg).
	\end{equation*}
	In particular, under Assumption~\ref{assump : np(n)}, we have for every $\eta > 0$
	\begin{equation*}
	\frac{1}{n} \log \mathbb{P} \bigg[ \frac{\vert\vert D^{(n)} \vert\vert_{\infty \rightarrow 1}}{n} > \eta \bigg] \longrightarrow -\infty, \qquad \text { as } n \rightarrow \infty.
	\end{equation*}
\end{lem}

\smallskip

\begin{lem} [Bernstein's inequality, Theorem~2.9 of \cite{BLM}]	\label{lem : Bernstein}
	Let $X_1, \cdots X_k$ be independent random variables with finite variance such that $X_i \le b$ for some $b > 0$ almost surely for each $i \in [k]$. Let $v = \sum_{i=1}^k \mathbb{E}[X_i^2]$, then we have
	\begin{equation*}
	\mathbb{P} \Big[ \sum_{i=1}^k \big( X_i - \mathbb{E}(X_i) \big) \ge u \Big]
	\le \exp \bigg( - \frac{u^2}{2(v + \frac{bu}{3})} \bigg).
	\end{equation*}
\end{lem}

\smallskip

When the interaction function $\phi$ belongs to a special class of functions, we shall see in the proof of Theorem~\ref{thm : W1 L} that the distance $W_1(L_n, \bar{L}_n)$ can be easily expressed in terms of the quantity $\vert\vert D^{(n)} \vert\vert_{\infty \rightarrow 1}$. This observation is inspired by the work of \cite{Oliveria:Reis}. To state more precisely, we introduce the notion of $L^1$-Fourier class of functions.

\smallskip

\begin{defn}	\label{def : Fourier class}
	Identifying $\mathbb{R}^{2d}$ with $\mathbb{R}^d \times \mathbb{R}^d$, we say that a function $f : \mathbb{R}^d \times \mathbb{R}^d \rightarrow \mathbb{R}$ belongs to the \textit{$L^1$-Fourier class}, if there exists a finite complex measure $m_{f}$ over $\mathbb{R}^{2d}$ such that for every $(x, y) \in \mathbb{R}^d \times \mathbb{R}^d$
	\begin{equation*}
	f(x, y) = \int_{\mathbb{R}^{2d}} \exp \Big(2\pi \sqrt{-1} \big\langle (x, y), \, z \big\rangle \Big) m_f(dz).
	\end{equation*}
\end{defn}

\smallskip

We recall that a finite complex measure $m$ over $\mathbb{R}^{2d}$ is a set function $m : \pazocal{B}(\mathbb{R}^{2d}) \rightarrow \mathbb{C}$ of the form $m = m_r^+ - m_r^- + \sqrt{-1}(m_i^+ - m_i^-)$, where each $m_r^+, m_r^-, m_i^+, m_i^-$ is a finite, $\sigma$-additive (nonnegative) measure over $\mathbb{R}^{2d}$. We define the total mass of $m$
\begin{equation*}
\vert\vert m \vert\vert_{TM} := m_r^+(\mathbb{R}^{2d}) + m_r^-(\mathbb{R}^{2d}) + m_i^+(\mathbb{R}^{2d}) + m_i^-(\mathbb{R}^{2d}).
\end{equation*}
If a function $f$ is an inverse $L^1$-transform of a function in $L^1(\mathbb{R}^{2d})$, then $f$ belongs to the $L^1$-Fourier class. In particular, any Schwartz function belongs to the $L^1$-Fourier class. An example of such function is the Kuramoto interaction; if $d=1$ and $\phi(x-y) = K \sin(y-x)$ for some constant $K$, then the corresponding complex measure is equal to
\begin{equation*}
m_{\phi} = \frac{K}{2\sqrt{-1}} \big( \delta_{(-1, 1)} + \delta_{(1, -1)} \big).
\end{equation*}
The finite system \eqref{eq : SDE} of ``oscillators'' with the Kuramoto interaction function is studied in \cite{Coppini}.

\bigskip

\section{Main results}
\label{sec: concentration}

This section consists of three parts. The first part shows that expectations of the $W_2$-distances between two empirical measures on $\mathbb{R}^d$ related to the systems \eqref{eq : Graphon} - \eqref{eq : SDE bar} converge to zero as the number of particles goes to infinity. The second part gives exponential bounds of the probabilities that Lipschitz function values of the particles $\bar{X}^n$ of the system \eqref{eq : SDE bar} on $\big(C([0, T] : \mathbb{R}^d)\big)^n$ deviate from their means; the stronger norm we use for the space $\big(C([0, T] : \mathbb{R}^d)\big)^n$, the stronger assumption on the initial distribution of the particles is needed. These results for the first two parts will be used in proving the results in the last subsection. In the last part, we derive several concentration results of the finite particle systems \eqref{eq : SDE}, \eqref{eq : SDE bar} toward the graphon particle system \eqref{eq : Graphon} under different metrics.

\medskip

\subsection{Concentration in mean of the $W_2$-distance}		\label{subsec : expectation}

Let us recall the law $\mu_{u, t}$ of \eqref{eq : Graphon}, the empirical measures \eqref{def : empirical measures at t} of the three systems and the averaged law $\widetilde{\mu}_t := \int_0^1 \mu_{u, t} \, du$ for every $t \in [0, T]$. We give two expectations converging to zero as $n \rightarrow \infty$ in the following. The proofs are provided in Section~\ref{subsec : expectation proof}.

\begin{prop}	\label{prop : expectation W2 L L bar}
	Under Assumptions~\ref{assump : Lipschitz} and \ref{assump : np(n)},
	\begin{equation*}
	\mathbb{E}\big[\sup_{0 \le t \le T} W_{2} (L_{n, t}, \bar{L}_{n, t}) \big] \longrightarrow 0,
	\end{equation*}
	as $n\rightarrow \infty$.
\end{prop}

\smallskip

\begin{prop}	\label{prop : expectation W2 L bar mu tilde}
	Under Assumptions~\ref{assump : Lipschitz}, \ref{assump : independent initial}, and \ref{assump : continuity}(c),
	\begin{equation*}
	\mathbb{E}\big[\sup_{0 \le t \le T} W_{2} (\widetilde{L}_{n, t}, \, \widetilde{\mu}_t) \big] \longrightarrow 0,
	\end{equation*}
	as $n\rightarrow \infty$.
\end{prop}

\smallskip

By virtue of Lemma~\ref{lem : LLN}, we have
\begin{align*}
\mathbb{E} \big[\sup_{0 \le t \le T} W^2_{2} (L_{n, t}, \widetilde{L}_{n, t}) \big]
&\le \mathbb{E} \Big[\sup_{0 \le t \le T} \frac{1}{n} \sum_{i=1}^n \big\vert X_{\frac{i}{n}}(t) - X^n_i(t) \big\vert^2 \Big]
\\
&\le \mathbb{E} \Big[ \frac{1}{n} \sum_{i=1}^n \big\Vert X_{\frac{i}{n}} - X^n_i \big\Vert_{\star, T}^2 \Big] 
\longrightarrow 0, \qquad \text{as } n \rightarrow \infty.
\end{align*}
Combining the last convergence with Propositions~\ref{prop : expectation W2 L L bar} and \ref{prop : expectation W2 L bar mu tilde}, we immediately have other convergences of the expectations.

\smallskip

\begin{cor}	\label{cor : expectations converging to zero}
	Under assumptions of Propositions~\ref{prop : expectation W2 L L bar}, \ref{prop : expectation W2 L bar mu tilde},
	\begin{equation*}
	\mathbb{E}\big[\sup_{0 \le t \le T} W_{2} (L_{n, t}, \widetilde{\mu}_t) \big] \longrightarrow 0,
	\qquad 
	\mathbb{E}\big[\sup_{0 \le t \le T} W_{2} (\bar{L}_{n, t}, \widetilde{\mu}_t) \big] \longrightarrow 0,
	\end{equation*}
	as $n\rightarrow \infty$.
\end{cor}

\medskip

\subsection{Concentration around mean}	\label{subsec : concentration around mean}

We present in this subsection the concentration of 1-Lipschitz function of the particles $\bar{X}^n$ around its mean, under two different norms $\ell^1$ and $\ell^2$. Proofs of the results rely on the transportation inequalities presented in Section~\ref{subsec : Transportation}, and they will be given in Section~\ref{subsec : concentration around mean proof}.

\smallskip

From Lemma~\ref{lem : Wasser equivalence}, we note that the condition~\eqref{con : indep initial with exponential} of Assumption~\ref{assump : independent initial} in Theorem~\ref{thm : concen ineq X bar1} below is equivalent to the condition
\begin{equation}	\label{con : initial wasser1}
	W_{1} (\mu_{u, 0}, \nu) \le \sqrt{2 \kappa H(\nu \vert \mu_{u, 0})},	\qquad \text{for every } u \in [0, 1], ~~~ \nu \in \pazocal{P}^1(S).
\end{equation}

\smallskip

\begin{thm}	\label{thm : concen ineq X bar1}
	Under Assumptions~\ref{assump : Lipschitz} and \ref{assump : independent initial}, there exists a constant $\delta > 0$, independent of $n$, such that for every $F \in Lip\Big(\big(C([0, T] : \mathbb{R}^d)\big)^n, \vert\vert\cdot\vert\vert_{n, 1}\Big)$ and every $a > 0$
	\begin{equation}	\label{eq : concentration ineq X bar1}
	\mathbb{P} \Big[ F(\bar{X}^n) - \mathbb{E}\big(F(\bar{X}^n)\big) > a \Big] \le 2 \exp \Big( -\frac{\delta a^2}{n} \Big)
	\end{equation}
	holds.
\end{thm}

\smallskip

We have the following analogous result to Theorem~\ref{thm : concen ineq X bar1}, when the condition \eqref{con : initial wasser1} is replaced by \eqref{con : initial wasser2}. For any $u \in [0, 1]$, if the initial law takes the form $\mu_{u, 0}(dx) = e^{-U(x)} dx$ for some $U \in C^2(\mathbb{R}^d)$ with Hessian bounded below in semidefinite order by $cI$ for some $c>0$, then $\mu_{u, 0}$ satisfies the condition \eqref{con : initial wasser2} with $\kappa = 1/c$. In particular, if $\mu_{u, 0}$ has the standard normal distribution on $\mathbb{R}^d$, then \eqref{con : initial wasser2} holds with $\kappa = 1$.
We emphasize that the concentration inequality \eqref{eq : concentration ineq X bar} is \textit{dimension-free}; the right-hand side does not depend on $n$.

\smallskip

\begin{thm}	\label{thm : concen ineq X bar2}
	Suppose that the initial particles $\{X_u(0)\}_{u \in [0, 1]}$ are independent with law $\mu_{u, 0} \in \pazocal{P}(\mathbb{R}^d)$, satisfying for some $\kappa > 0$
	\begin{equation}	\label{con : initial wasser2}
		W_{2} (\mu_{u, 0}, \nu) \le \sqrt{2 \kappa H(\nu \vert \mu_{u, 0})},	\qquad \text{for every } u \in [0, 1], ~~~ \nu \in \pazocal{P}^2(S).
	\end{equation}
	Under Assumption~\ref{assump : Lipschitz}, there exists a constant $\delta > 0$, independent of $n$, such that for every $F \in Lip\Big(\big(C([0, T] : \mathbb{R}^d)\big)^n, \vert\vert\cdot\vert\vert_{n, 2}\Big)$ and every $a > 0$
	\begin{equation}	\label{eq : concentration ineq X bar}
	\mathbb{P} \Big[ F(\bar{X}^n) - \mathbb{E}\big(F(\bar{X}^n)\big) > a \Big] \le 2 \exp ( -\delta a^2)
	\end{equation}
	holds.
\end{thm}

\medskip

\subsection{Concentration toward the graphon system}	\label{subsec : concentration toward GS}

Recalling the notations in \eqref{def : empirical measures at t}, we now provide concentration in terms of (1 and 2)-Wasserstein distance of the empirical measures of the finite particle systems toward the averaged measure $\widetilde{\mu}_t$ of the graphon system. Proofs will be given in Section~\ref{subsec : concentration toward GS proof}

\smallskip

First, we have the following concentration result of $\bar{L}_{n, t}$ toward $\widetilde{\mu}_t$ in terms of the $W_1$-distance, due to Theorem~\ref{thm : concen ineq X bar1}. 

\smallskip

\begin{thm}		\label{thm : W1 L bar}
	Under Assumptions~\ref{assump : Lipschitz}, \ref{assump : independent initial}, and \ref{assump : continuity}(c), there exist constants $\delta > 0$, which is independent of $n$, and $N \in \mathbb{N}$ such that
	\begin{equation}	\label{ineq : Wasserstein 2 bound for bar L1}
	\mathbb{P} \Big[ \sup_{0 \le t \le T} W_{1} (\bar{L}_{n, t}, \widetilde{\mu}_t) > a \Big] \le 2 \exp \Big( - \frac{\delta a^2n}{4} \Big)
	\end{equation}
	holds for every $a > 0$ and every $n \ge N$.
\end{thm}

\smallskip

\begin{rem}
	Theorem~\ref{thm : W1 L bar} gives the same exponential bound as in Theorem~2.1 of \cite{BW}. Their proof mainly focuses on computing certain sub-Gaussian estimates, whereas our argument relies on the concentration property \eqref{eq : concentration ineq X bar1} of the system \eqref{eq : SDE bar}. Applying the same latter argument, we can even deduce the exponential bound in terms of the $W_2$-metric in Theorem~\ref{thm : W2 L bar} below.
\end{rem}

\smallskip

In Section~\ref{subsec: EE finite}, the concept of $L^1$-Fourier class, along with Bernstein's inequality was introduced to express $W_1(L_n, \bar{L}_n)$ in terms of $\Vert D^{(n)} \Vert_{\infty \rightarrow 1}$. This gives rise to the following concentration result of the particle system \eqref{eq : SDE} toward the graphon system.

\smallskip

\begin{thm}		\label{thm : W1 L}
	Suppose that the components of the interaction function $\phi$ belong to the $L^1$-Fourier class~(Definition~\ref{def : Fourier class}). Under Assumptions~\ref{assump : Lipschitz}, \ref{assump : independent initial}, \ref{assump : continuity}(c), and \ref{assump : np(n)}, there exist constants $\delta > 0$, which is independent of $n$, and $N \in \mathbb{N}$ such that
	\begin{equation}	\label{ineq : Wasserstein 1 bound for L1}
	\mathbb{P} \Big[ \sup_{0 \le t \le T} W_{1} (L_{n, t}, \widetilde{\mu}_t) > a \Big]
	\le 3 \exp \Big( - \frac{\delta a^2n}{16} \Big)
	\end{equation}
	holds for every $a > 0$ and every $n \ge N$. For general interaction functions $\phi$, we have instead
	\begin{equation}	\label{ineq : Wasserstein 1 bound for BL}
	\mathbb{P} \Big[ \sup_{0 \le t \le T} d_{BL} (L_{n, t}, \widetilde{\mu}_t) > a \Big]
	\le 3 \exp \Big( - \frac{\delta a^2n}{16} \Big).
	\end{equation}
\end{thm}

\smallskip

The following result gives the concentration of $\bar{L}_{n, t}$ toward $\widetilde{\mu}_t$ as Theorem~\ref{thm : W1 L bar}, but in terms of the $W_2$-metric. Its proof is similar to that of Theorem~\ref{thm : W1 L bar}, but Theorem~\ref{thm : concen ineq X bar2} is used in place of Theorem~\ref{thm : concen ineq X bar1}.

\smallskip

\begin{thm}		\label{thm : W2 L bar}
	Under Assumptions~\ref{assump : Lipschitz}, \ref{assump : independent initial}, and \ref{assump : continuity}, together with the condition \eqref{con : initial wasser2}, there exist constants $\delta > 0$, independent of $n$, and $N \in \mathbb{N}$ such that
	\begin{equation}	\label{ineq : Wasserstein 2 bound for bar L}
		\mathbb{P} \big[ \sup_{0 \le t \le T} W_2 (\bar{L}_{n, t}, \widetilde{\mu}_t) > a \big] \le 2 \exp \Big( - \frac{\delta a^2n}{4} \Big)
	\end{equation}
	holds for every $a > 0$ and every $n \ge N$.
\end{thm}

\smallskip

Since we have the exponential bound in \eqref{ineq : Wasserstein 2 bound for bar L} in the $W_2$-metric, one naturally expects to obtain a similar bound to \eqref{ineq : Wasserstein 1 bound for L1} in the $W_2$-metric as well. In order to achieve this, we need to find the exponential bound for the probability $\mathbb{P} \big[ \sup_{0 \le t \le T} W_2 (L_{n, t}, \bar{L}_{n, t}) > a \big]$, which requires us to handle the quantity $\Vert (D^{(n)})^\top D^{(n)} \Vert_{\infty \rightarrow 1}$, instead of $\Vert D^{(n)} \Vert_{\infty \rightarrow 1}$ as in the proof of Theorem~\ref{thm : W1 L}. Controlling this quantity is done in Lemma~\ref{lem : DTD} under an extra condition on the sparsity parameter $p(n)$, more restrictive condition than the one in Assumption~\ref{assump : np(n)}.

\smallskip

\begin{assump}	\label{assump : np(n)^2}
	The sparsity parameter sequence $\{p(n)\}_{n \in \mathbb{N}} \subset (0, 1]$ of the system \eqref{eq : SDE} satisfies either one of the following:
	\begin{enumerate} [(a)]
		\item $p(n) \rightarrow 0$ and $np(n)^2 \rightarrow \infty$ as $n \rightarrow \infty$, or
		\item $p(n) \equiv 1$ for every $n \in \mathbb{N}$.
	\end{enumerate}
\end{assump}

\smallskip

Recalling the notations of \eqref{def : P bar}, the following lemma is needed when proving Theorem~\ref{thm : W2 L}. Its proof in Section~\ref{subsec : concentration toward GS proof} is similar to that of Lemma~\ref{lem : cut norm}, but requires more involved applications of Bernstein's inequality.

\begin{lem}	\label{lem : DTD}
	Under Assumption~\ref{assump : np(n)^2}, there exists $N \in \mathbb{N}$ such that 
	\begin{equation*}
	\mathbb{P} \bigg[ \frac{\Vert (D^{(n)})^\top D^{(n)} \Vert_{\infty \rightarrow 1}}{n} > \eta \bigg] 
	\le 3n^2 \exp \bigg( - \frac{2\eta^2 np(n)^4}{9+4\eta} \bigg)
	\end{equation*}
	holds for every $n \ge N$ and $\eta > 0$.
\end{lem}

\smallskip

\begin{thm}		\label{thm : W2 L}
	Suppose that the components of the interaction function $\phi$ belong to the $L^1$-Fourier class. Under Assumptions~\ref{assump : Lipschitz}, \ref{assump : independent initial}, \ref{assump : continuity}, and \ref{assump : np(n)^2}(a), together with the condition \eqref{con : initial wasser2}, there exist constants $K > 0$, which is independent of $n$, and $N \in \mathbb{N}$ such that
	\begin{equation}	\label{ineq : Wasserstein 2 bound for L1}
	\mathbb{P} \Big[ \sup_{0 \le t \le T} W_{2} (L_{n, t}, \widetilde{\mu}_t) > a \Big]
	\le 4n^2 \exp \bigg( - \frac{a^4 np(n)^4}{72K^2+8a^2K} \bigg)
	\end{equation}
	holds for every $a > 0$ and every $n \ge N$.
	 
	\smallskip
	 
	Furthermore, if Assumption~\ref{assump : np(n)^2}(a) is replaced by Assumption~\ref{assump : np(n)^2}(b), we have the exponential bound in $n$; there exist constants $\delta > 0$, which is independent of $n$, and $N \in \mathbb{N}$ such that
	\begin{equation}	\label{ineq : Wasserstein 2 bound for L1 exponential}
	\mathbb{P} \Big[ \sup_{0 \le t \le T} W_{2} (L_{n, t}, \widetilde{\mu}_t) > a \Big]
	\le \exp \bigg( - \frac{\delta a^4 n}{a^2+ \delta} \bigg)
	\end{equation}
	holds for every $a > 0$ and every $n \ge N$.
\end{thm}

\bigskip

\section{Proofs} 	\label{sec : proofs}

In this section, we provide proofs of the results stated in Section~\ref{sec: concentration}.

\medskip
 
\subsection{Proofs of results in Section~\ref{subsec : expectation}}	\label{subsec : expectation proof}

\subsubsection{Proof of Proposition~\ref{prop : expectation W2 L L bar}}
	Let us recall the identity \eqref{eq : diff X X bar}, along with the notations \eqref{def : P bar}. Using H\"older's inequality, there exists $K > 0$, depending on $\phi$ and $\psi$, such that for every $t \in [0, T]$
	\begin{align*}
	\sup_{0 \le t \le T} \big\vert X^n_i(t) - \bar{X}^n_i(t) \big\vert^2
	&\le KT \int_0^T \Big\vert \sum_{j=1}^n D^{(n)}_{i, j} \phi\big( X^n_i(s), X^n_j(s) \big) \Big\vert^2 \, ds
	\\
	& + KT \int_0^T \Big\vert \sum_{j=1}^n \bar{P}^{(n)}_{i, j} \Big ( \phi\big(X^n_i(s), X^n_j(s) \big) - \phi\big(\bar{X}^n_i(s), \bar{X}^n_j(s) \big) \Big) \Big\vert^2 \, ds
	\\
	&+ KT \int_0^T \big\vert X^n_i(s) - \bar{X}^n_i(s) \big\vert^2 \, ds.
	\end{align*}
	Taking the expectation to the first term, and using the independence of $\{D^{(n)}_{i, j}\}_{j \in [n]}$ and the boundedness of $\phi$, we have
	\begin{align*}
	KT \, \mathbb{E} \int_0^T \Big\vert \sum_{j=1}^n D^{(n)}_{i, j} \phi\big( X^n_i(s), X^n_j(s) \big) \Big\vert^2 \, ds
	\le KT \int_0^T \sum_{j=1}^n \mathbb{E}[(D^{(n)}_{i, j})^2] \, ds
	\le \frac{KT^2}{np(n)}.
	\end{align*}
	For the second term, H\"older's inequality and the Lipschitz continuity of $\phi$ give
	\begin{align*}
	KT \, \mathbb{E} \int_0^T \Big\vert \sum_{j=1}^n \bar{P}^{(n)}_{i, j} &\Big ( \phi\big(X^n_i(s), X^n_j(s) \big) - \phi\big(\bar{X}^n_i(s), \bar{X}^n_j(s) \big) \Big) \Big\vert^2 \, ds
	\\
	&\le KT \, \mathbb{E} \int_0^T \frac{1}{n}\sum_{j=1}^n \big( \vert X^n_i(s) - \bar{X}^n_i(s) \vert^2 + \vert X^n_j(s) - \bar{X}^n_j(s) \vert^2 \big) \, ds.
	\end{align*}
	Combining above inequalities and averaging over $i \in [n]$, we obtain
	\begin{align*}
	\frac{1}{n} \sum_{i=1}^n \mathbb{E} \big[ \sup_{0 \le t \le T} \vert X^n_i(t) - \bar{X}^n_i(t) \vert^2 \big]
	\le \frac{KT^2}{np(n)} + KT\int_0^T \frac{1}{n} \sum_{i=1}^n \mathbb{E} \big[ \sup_{0 \le u \le s}\vert X^n_i(u) - \bar{X}^n_i(u) \vert^2 \big] \, ds.
	\end{align*}
	Gr\"onwall's inequality yields
	\begin{equation*}
	\frac{1}{n} \sum_{i=1}^n \mathbb{E} \big[ \sup_{0 \le t \le T} \vert X^n_i(t) - \bar{X}^n_i(t) \vert^2 \big]
	\le \frac{KT^2 \exp(KT^2)}{np(n)},
	\end{equation*}
	and thus
	\begin{align*}
	\mathbb{E}\big[\sup_{0 \le t \le T} W^2_{2} (L_{n, t}, \bar{L}_{n, t}) \big]
	&\le \mathbb{E} \big[ \sup_{0 \le t \le T} \frac{1}{n}\sum_{i=1}^n \vert X^n_i(t) - \bar{X}^n_i(t) \vert^2 \big]
	\\
	&\le \frac{1}{n} \sum_{i=1}^n \mathbb{E} \big[ \sup_{0 \le t \le T} \vert X^n_i(t) - \bar{X}^n_i(t) \vert^2 \big]
	\le \frac{KT^2 \exp(KT^2)}{np(n)} \longrightarrow 0, \text{ as } n \rightarrow \infty.
	\end{align*}
\hfill $\qed$

\smallskip

\subsubsection{Proof of Proposition~\ref{prop : expectation W2 L bar mu tilde}}
	We partition the interval $[0, T]$ into $M := \lceil \frac{T}{\Delta} \rceil$ subintervals of length $\Delta > 0$:
	\begin{equation*}
	[0, T] = [0, \Delta] \cup [\Delta, 2\Delta] \cup \cdots \cup [(M-1)\Delta, T] 
	=: \cup_{h=1}^{M} \Delta_h,
	\end{equation*}
	where $\Delta_h := [(h-1)\Delta, \, h\Delta]$ for $h = 1, \cdots, M-1$ and $\Delta_M = [(M-1)\Delta, \, T]$, and we choose the value of $\Delta$ later. With the notation
	\begin{equation*}
	\widetilde{\mu}_{n, t} := \frac{1}{n} \sum_{i=1}^n \mu_{\frac{i}{n}, t} = \frac{1}{n} \sum_{i=1}^n \pazocal{L}\big(X_{\frac{i}{n}}(t)\big),
	\end{equation*}
	triangle inequality gives
	\begin{align*}
	\mathbb{E} \big[ \sup_{0 \le t \le T} W_{2} (\widetilde{L}_{n, t}, \, \widetilde{\mu}_t) \big]
	&= \mathbb{E} \Big[ \max_{h \in [M]} \, \sup_{t \in \Delta_h} W_{2} (\widetilde{L}_{n, t}, \, \widetilde{\mu}_t) \Big]
	\\
	& \le \mathbb{E} \Big[ \max_{h \in [M]} \, \sup_{t \in \Delta_h} W_{2} (\widetilde{L}_{n, t}, \, \widetilde{L}_{n, (h-1)\Delta}) \Big]
	+ \mathbb{E} \Big[ \sup_{h \in [M]} W_{2} (\widetilde{L}_{n, (h-1)\Delta}, \, \widetilde{\mu}_{n, (h-1)\Delta}) \Big]
	\\
	& \qquad + \mathbb{E} \Big[ \sup_{h \in [M]} W_{2} (\widetilde{\mu}_{n, (h-1)\Delta}, \, \widetilde{\mu}_{(h-1)\Delta}) \Big]
	+ \mathbb{E} \Big[ \sup_{h \in [M]} \, \sup_{t \in \Delta_h} W_{2} (\widetilde{\mu}_{(h-1)\Delta}, \, \widetilde{\mu}_{t}) \Big].
	\\
	& =: E_1 + E_2 + E_3 + E_4.
	\end{align*}
	For the first term $E_1$, we note that there exists $K>0$, depending on the bounds of $\phi$, $\psi$, and $\sigma$, such that 
	\begin{equation*}
	\big\vert X_{\frac{i}{n}}(s) - X_{\frac{i}{n}}(u) \big\vert^2 \le K \vert s-u \vert^2 + K \big\vert B_{\frac{i}{n}}(s) - B_{\frac{i}{n}}(u) \big\vert^2
	\end{equation*}
	holds for every $0 \le u \le s \le T$, and thus we have
	\begin{align*}
	(E_1)^2 
	&\le \mathbb{E} \Big[ \max_{h \in [M]} \, \sup_{t \in \Delta_h} W^2_{2} (\widetilde{L}_{n, t}, \, \widetilde{L}_{n, (h-1)\Delta}) \Big]
	\le	\mathbb{E} \Big[ \max_{h \in [M]} \, \sup_{t \in \Delta_h} \frac{1}{n}\sum_{i=1}^n \big\vert X_{\frac{i}{n}}(t) - X_{\frac{i}{n}}\big((h-1)\Delta\big) \big\vert^2 \Big]
	\\
	& \le \mathbb{E} \bigg[ \max_{h \in [M]} \, \sup_{t \in \Delta_h} \frac{1}{n}\sum_{i=1}^n \bigg( K \Delta^2 + K \Big\vert B_{\frac{i}{n}}(t) - B_{\frac{i}{n}}\big((h-1)\Delta\big) \Big\vert^2 \bigg) \bigg]
	\\
	& \le K \Delta^2 + K\mathbb{E} \Big[ \max_{h \in [M]} \, \sup_{t \in \Delta_h} \frac{1}{n}\sum_{i=1}^n \big\vert B_{\frac{i}{n}}(t) - B_{\frac{i}{n}}\big((h-1)\Delta\big) \big\vert^2 \Big].
	\end{align*}
	Applying H\"older's inequality twice, the last expectation is bounded above by
	\begin{align*}
	& \quad \mathbb{E} \bigg[ \max_{h \in [M]} \, \sup_{t \in \Delta_h} \Big( \frac{1}{n} \sum_{i=1}^n \big\vert B_{\frac{i}{n}}(t) - B_{\frac{i}{n}}\big((h-1)\Delta\big) \big\vert^4 \Big)^{\frac{1}{2}} \bigg]
	\\
	&\le \sqrt{\mathbb{E} \Big[ \max_{h \in [M]} \, \sup_{t \in \Delta_h} \frac{1}{n}\sum_{i=1}^n \big\vert B_{\frac{i}{n}}(t) - B_{\frac{i}{n}}\big((h-1)\Delta\big) \big\vert^4 \Big]}
	\\
	&\le \sqrt{ \sum_{h \in [M]} \frac{1}{n}\sum_{i=1}^n \mathbb{E} \Big[ \sup_{t \in \Delta_h}  \big\vert B_{\frac{i}{n}}(t) - B_{\frac{i}{n}}\big((h-1)\Delta\big) \big\vert^4 \Big]}
	\le \sqrt{M C_4 \mathbb{E}[\Delta^2]} \le \sqrt{(T+1) C_4 \Delta}. 
	\end{align*}
	The second-last inequality uses the properties of the increments of Brownian motion and the Burkholder-Davis-Gundy inequality with the positive constant $C_4$. Therefore, we have the bound
	\begin{equation*}
	(E_1)^2 \le K \Delta^2 + K \sqrt{(T+1)C_4 \Delta}.
	\end{equation*}
	
	\smallskip
	
	For the second expectation $E_2$, H\"older's inequality and Lemma~\ref{lem : independent empirical measures} give
	\begin{align*}
	E_2 
	&\le \sum_{h \in [M]} \mathbb{E} \big[ W_{2} (\widetilde{L}_{n, (h-1)\Delta}, \, \widetilde{\mu}_{n, (h-1)\Delta}) \big]
	\le \sum_{h \in [M]} \Big( \mathbb{E} \big[ W^3_{3} (\widetilde{L}_{n, (h-1)\Delta}, \, \widetilde{\mu}_{n, (h-1)\Delta}) \big] \Big)^{1/3}
	\\
	&\le K \sum_{h \in [M]} \bigg( \int_{\mathbb{R}^d} \vert x \vert^4 \widetilde{\mu}_{n, (h-1)\Delta}(dx) \bigg)^{1/4} \alpha_{3, 4}^{1/3}(n)
	\le 
	KM \alpha_{3, 4}^{1/3}(n) \longrightarrow 0,
	\end{align*}
	as $n \rightarrow \infty$, where the last inequality follows from Lemma~\ref{lem : fourth moment}.
	
	\smallskip
	
	On the other hand, using the convexity of $W^2_2(\, \cdot \, , \, \cdot \, )$ and Lemma~\ref{lem : continuity of graphon system}(b), there exists $K > 0$ satisfying
	\begin{equation*}
	E_3^2 \le \mathbb{E} \Big[ \sup_{h \in [M]} W^2_{2} (\widetilde{\mu}_{n, (h-1)\Delta}, \, \widetilde{\mu}_{(h-1)\Delta}) \Big] 
	\le \mathbb{E} \Big[ \sup_{h \in [M]} \int_0^1 W^2_{2} (\mu_{\frac{\lceil nu \rceil}{n}, (h-1)\Delta}, \, \mu_{u, (h-1)\Delta}) \, du \Big]
	\le \frac{K}{n^2}.
	\end{equation*}
	
	\smallskip
	
	Finally, for the last term $E_4$, we note from a straightforward computation that there exists $K > 0$ satisfying
	\begin{equation*}
	\mathbb{E} \big\vert X_u(t) - X_u(s) \big\vert^2 \le K \vert t-s \vert^2 + K\mathbb{E} \big\vert B_u(t) - B_u(s) \big\vert^2 \le K \vert t-s \vert
	\end{equation*}
	for every $u \in [0, 1]$ and $s, t \in [0, T]$ satisfying $\vert t-s \vert \le 1$. Thus, we have
	\begin{align*}
	E^2_4 &\le \sup_{h \in [M]} \, \sup_{t \in \Delta_h} \int_0^1 W^2_{2} (\mu_{u, (h-1)\Delta}, \, \mu_{u, t}) \, du
	\\
	&\le \sup_{h \in [M]} \, \sup_{t \in \Delta_h} \int_0^1 \mathbb{E} \big\vert X_u \big( (h-1)\Delta \big) - X_u(t) \big\vert^2 \, du
	\le K \Delta.
	\end{align*}
	
	\smallskip
	
	Let us combine all the bounds from $E_1$ to $E_4$. For any given $\epsilon > 0$, we can choose small enough $\Delta$ such that $E_1 + E_4 < \epsilon/2$. Then, we can choose large enough $N \in \mathbb{N}$ satisfying $E_2 + E_3 < \epsilon/2$ for every $n \ge N$, which implies $\mathbb{E} \big[ \sup_{0 \le t \le T} W_{2} (\widetilde{L}_{n, t}, \, \widetilde{\mu}_t) \big] < \epsilon$ for every $n \ge N$.
\hfill $\qed$

\medskip

\subsection{Proofs of results in Section~\ref{subsec : concentration around mean}}	\label{subsec : concentration around mean proof}

\subsubsection{Proof of Theorem~\ref{thm : concen ineq X bar1}}
	Let us fix an arbitrary $n \in \mathbb{N}$. We shall naturally identify elements of $(\mathbb{R}^d)^n$ with those of $\mathbb{R}^{dn}$, identify elements of $\big( C([0, T] : \mathbb{R}^d)\big)^n$ with those of $C\big([0, T] : (\mathbb{R}^d)^n\big)$, and we shall specify which norm we use for each space. We can express the SDE \eqref{eq : SDE bar} in the form of \eqref{eq : SDE form} with $k = dn$, by setting
	\begin{enumerate} [(1)]
		\item $(\mathbb{R}^d)^n \ni x = (x_i)_{i \in [n]}$, where $x_i = X_{i/n}(0)$;
		\item $C\big([0, T] : (\mathbb{R}^d)^n\big) \ni X^x = (\bar{X}^n_i)_{i \in [n]}$ where $\bar{X}^n_i = (\bar{X}^n_{i, k})_{k \in [d]}$;
		\item $b : [0, T] \times C\big([0, T] : (\mathbb{R}^d)^n\big) \rightarrow (\mathbb{R}^d)^n$ such that $b = (b_i)_{i \in [n]}$, $b_i = (b_{i, k})_{k \in [d]}$ where \\
		$b_{i, k}(t, X^x) = \frac{1}{n} \sum_{j=1}^n \phi_k\big(\bar{X}^n_i(t), \bar{X}^n_j(t)\big)G(\frac{i}{n}, \frac{j}{n}) + \psi_k\big(\bar{X}^n_i(t)\big)$;
		\item $W = (W_i)_{i \in [n]}$ is a $(dn)$-dimensional Brownian motion, where $W_i \equiv B_{i/n}$;
		\item $\Sigma$ is a block-diagonal $(dn) \times (dn)$ matrix with block diagonal entries $\sigma$.
	\end{enumerate}
	In order to apply Lemma~\ref{lem : wasser entropy}, it suffices to check the condition \eqref{con : b lipschitz}; for any $X, Y \in C([0, T] : \mathbb{R}^{dn})$, H\"older inequality and the Lipschitz continuity of $\phi$, $\psi$ indeed yield for every $t \in [0, T]$
	\begin{align*}
	& \qquad \big\vert b(t, X) - b(t, Y) \big\vert^2
	\\
	&= \sum_{i=1}^n \sum_{k=1}^d \bigg \vert \frac{1}{n} \sum_{j=1}^n \Big(\phi_k \big(X_i(t), X_j(t)\big) - \phi_k \big(Y_i(t), Y_j(t)\big) \Big) G(\frac{i}{n}, \frac{j}{n}) + \psi_k\big(X_i(t) \big) - \psi_k\big(Y_i(t)\big) \bigg\vert^2
	\\
	& \le 2 \sum_{i=1}^n \sum_{k=1}^d \bigg[ \Big \vert \frac{1}{n} \sum_{j=1}^n \Big(\phi_k \big(X_i(t), X_j(t)\big) - \phi_k \big(Y_i(t), Y_j(t)\big) \Big) G(\frac{i}{n}, \frac{j}{n}) \Big\vert^2 + \Big\vert \psi_k\big(X_i(t) \big) - \psi_k\big(Y_i(t)\big) \Big \vert^2 \bigg]
	\\
	&\le 2 \sum_{i=1}^n \bigg[ K^2  \big\vert X_i(t) - Y_i(t) \big\vert^2 + \sum_{k=1}^d \frac{1}{n} \sum_{j=1}^n \Big\vert \phi_k \big(X_i(t), X_j(t)\big) - \phi_k \big(Y_i(t), Y_j(t)\big) \Big\vert^2 \bigg]
	\\
	&\le (4K^2+1) \sum_{i=1}^n \big\vert X_i(t) - Y_i(t) \big\vert^2
	= (4K^2+1) \sum_{i=1}^n \sum_{k=1}^d \big\vert X_{i, k}(t) - Y_{i, k}(t) \big\vert^2
	\le (4K^2+1) \vert\vert X-Y \vert\vert^2_{\star, t}.
	\end{align*}
	Let $P^x \in \pazocal{P}\big(C([0, T] : \mathbb{R}^{dn})\big)$ be the law of the solution of \eqref{eq : SDE bar} in the notations of (1)-(5) above, then we have from Lemma~\ref{lem : wasser entropy} for any $Q \in \pazocal{P}\big(C([0, T] : \mathbb{R}^{dn})\big)$, 
	\begin{equation*}
	W_{1, \big(C([0, T] : \mathbb{R}^{dn}), \, \vert\vert\cdot\vert\vert_{dn, 2}\big)} (P^x, Q) 
	\le \sqrt{2c_1 H(Q \vert P^x)},
	\end{equation*}
	for some $c_1 > 0$.
	
	\smallskip
	
	For an arbitrary $F \in Lip\Big(\big(C([0, T] : \mathbb{R}^d)\big)^n, \vert\vert\cdot\vert\vert_{n, 1}\Big)$, H\"older's inequality shows that $F$ is $\sqrt{n}$-Lipschitz function of the space $\big(C([0, T] : \mathbb{R}^{dn}), \, \Vert \cdot \Vert_{dn, 2}\big)$; we indeed obtain for $X, Y \in \Big(\big(C([0, T] : \mathbb{R}^d)\big)^n, \Vert \cdot \Vert_{n, 1}\Big)$
	\begin{align*}
	\big\vert F(X) - F(Y) \big\vert 
	&\le \Vert X-Y \Vert_{n, 1} 
	= \sum_{i=1}^d \Vert X_i - Y_i \Vert_{\star, T}
	\\
	&\le \sqrt{n \sum_{i=1}^d \Vert X_i - Y_i\Vert_{\star, T}^2}
	\le \sqrt{n \sum_{i=1}^d \sum_{k=1}^d \Vert X_{i, k} - Y_{i, k}\Vert_{\star, T}^2}
	= \sqrt{n} \, \Vert X-Y \Vert_{dn, 2}.
	\end{align*}
	Thus, Lemma~\ref{lem : Wasser equivalence} implies
	\begin{equation}	\label{eq : P^x F bound1}
	P^x \big(F - \langle P^x, F \rangle > a \big) \le \exp\Big(-\frac{a^2}{2c_1 n}\Big),
	\end{equation}
	for any $a > 0$.
	
	\smallskip
	
	We now claim that there exists a positive constant $c_2$, which does not depend on $n$, such that the map $x \mapsto \langle P^x, F \rangle$ is $c_2$-Lipschitz on $(\mathbb{R}^d)^n$ with respect to the Euclidean $\ell^p$-norm for any $F \in Lip\Big(\big(C([0, T] : \mathbb{R}^d)\big)^n, \vert\vert\cdot\vert\vert_{n, p}\Big)$ and for any $p = 1, 2$. Given any $x, y \in (\mathbb{R}^d)^n$, we couple $P^x$ and $P^y$ by solving the system \eqref{eq : SDE bar} from the two initial states $x, y$ with the same Brownian motion, and denote the coupling by $\pi_{x, y}$. We deduce that for $\pazocal{L}(X) = P^x$ and $\pazocal{L}(Y) = P^y$
	\begin{align}	\label{ineq : coupling x y}
	\big\vert \langle P^x, F \rangle - \langle P^y, F \rangle \big\vert^p
	\le \int \big\vert F(X) - F(Y) \big\vert^p \pi_{x, y}(dX, dY)
	\le \int \Vert X - Y \Vert^p_{n, p} \, \pi_{x, y}(dX, dY).
	\end{align}
	When $p = 2$, we use a standard argument~(the trivial inequality $(a+b)^2 \le 2(a^2+b^2)$, the Lipschitz continuity from Assumption~\ref{assump : Lipschitz}(a), and a series of H\"older inequalities) to derive
	\begin{align*}
	\sum_{i=1}^n \vert X_i(t) - Y_i(t) \vert^2
	&\le 2 \vert\vert x-y\vert\vert^2_{n, 2} + 2Kt \sum_{i=1}^n \int_0^t \Big( \vert X_i(s) - Y_i(s) \vert + \frac{1}{n}\sum_{j=1}^n \vert X_j(s) - Y_j(s) \vert \Big)^2 \, ds
	\\
	&\le 2 \vert\vert x-y\vert\vert^2_{n, 2} + 8Kt \int_0^t \sum_{i=1}^n \vert X_i(s) - Y_i(s) \vert^2 \, ds, \qquad \forall \, t \in [0, T].
	\end{align*}
	Gr\"onwall's inequality yields that the last integrand of \eqref{ineq : coupling x y} for $p=2$ is bounded by
	\begin{equation*}
	\vert\vert X-Y \vert\vert^2_{n, 2} \le c_2^2 \vert\vert x-y \vert\vert^2_{n, 2}
	\end{equation*}
	for some constant $c_2 > 0$, which depends on $\phi$, $\psi$, and $T$, but not on $n$. When $p=1$, proving $\vert\vert X-Y \vert\vert_{n, 1} \le c_2 \vert\vert x-y \vert\vert_{n, 1}$ is easier, and the claim follows.
	
	\smallskip
	
	On the other hand, we apply Lemmas~\ref{lem : Wasser tensor}, \ref{lem : Wasser equivalence} to the assumption \eqref{con : initial wasser1} to obtain for every $f \in Lip\big( (\mathbb{R}^d)^n, \, \vert\vert\cdot\vert\vert_{n, 1} \big)$ and for any $a > 0$
	\begin{equation}	\label{eq : mu_0 bound1}
	\mu^n_0 \big(f - \langle \mu^n_0, f \rangle > a \big) \le \exp\Big(-\frac{a^2}{2\kappa n}\Big).
	\end{equation}
	
	\smallskip
	
	We conclude from \eqref{eq : P^x F bound1}, the above claim, and \eqref{eq : mu_0 bound1}
	\begin{align*}
	\mathbb{P} \Big[ F(\bar{X}^n) - \mathbb{E}\big(F(\bar{X}^n)\big) > a \Big]
	& \le \mathbb{E} \Big[ \mathbb{P} \big( F(\bar{X}^n) - \langle P^{x}, F \rangle > \frac{a}{2} \, \big \vert \, \bar{X}^n(0) = x \big) \Big]
	\\
	& \qquad \qquad \qquad + \mathbb{P} \Big( \langle P^{\bar{X}^n(0)}, F \rangle - \mathbb{E} \big[\langle P^{\bar{X}^n(0)}, F \rangle \big] > \frac{a}{2} \Big)
	\\
	& \le \exp\Big(-\frac{a^2}{8c_1 n}\Big) + \exp\Big(-\frac{a^2}{8\kappa c_2^2n}\Big).
	\end{align*}
	The assertion \eqref{eq : concentration ineq X bar} follows by choosing $1 / \delta = 8 \max(c_1, \kappa c_2^2)$.
\hfill $\qed$

\smallskip

\subsubsection{Proof of Theorem~\ref{thm : concen ineq X bar2}}
	We follow the proof of Theorem~\ref{thm : concen ineq X bar1}. Identifying the elements of $\big(C([0, T] : \mathbb{R}^d)\big)^n$ with those of $C([0, T] : \mathbb{R}^{dn})$, expressing the SDE \eqref{eq : SDE bar} in the form of \eqref{eq : SDE form}, and applying Lemma~\ref{lem : wasser entropy}, there exists a positive constant $c_1 > 0$ such that
	\begin{equation*}
	W_{1, (C([0, T] : \mathbb{R}^{dn}), \, \vert\vert\cdot\vert\vert_{dn, 2})}(P^x, Q) \le \sqrt{2c_1 H(Q \vert P^x)}
	\end{equation*}
	holds for any $Q \in \pazocal{P}\big(C([0, T] : \mathbb{R}^{dn})\big)$. Here, $P^x$ is the law of the solution of \eqref{eq : SDE bar}. Moreover, Lemma~\ref{lem : Wasser equivalence} implies
	\begin{equation}	\label{eq : P^x F bound}
	P^x \big(F - \langle P^x, F \rangle > a \big) \le \exp\Big(-\frac{a^2}{2c_1}\Big),
	\end{equation}
	for any $a > 0$ and every $F \in Lip\big(C([0, T] : \mathbb{R}^{dn}), \Vert\cdot\Vert_{dn, 2}\big)$. It is easy to check that every function in $Lip\big(C([0, T] : \mathbb{R}^{dn}), \Vert\cdot\Vert_{dn, 2}\big)$ also belongs to $Lip\Big(\big(C([0, T] : \mathbb{R}^d)\big)^n, \Vert\cdot\Vert_{n, 2}\Big)$, thus the inequality \eqref{eq : P^x F bound} also holds for every $F \in Lip\Big(\big(C([0, T] : \mathbb{R}^d)\big)^n, \Vert\cdot\Vert_{n, 2}\Big)$.
	
	\smallskip
	
	We now apply Lemmas~\ref{lem : Wasser tensor}(ii), \ref{lem : Wasser equivalence} to the assumption \eqref{con : initial wasser2} to deduce
	\begin{equation}	\label{eq : mu_0 bound}
	\mu^n_0 \big(f - \langle \mu^n_0, f \rangle > a \big) \le \exp\Big(-\frac{a^2}{2\kappa}\Big),
	\end{equation}
	for every $f \in Lip\big( (\mathbb{R}^d)^n, \, \vert\vert\cdot\vert\vert_{n, 2} \big)$ and for any $a > 0$.
	
	\smallskip
	
	From \eqref{eq : P^x F bound}, \eqref{eq : mu_0 bound}, and the claim in the proof of Theorem~\ref{thm : concen ineq X bar1}, we conclude that
	\begin{align*}
	\mathbb{P} \Big[ F(\bar{X}^n) - \mathbb{E}\big(F(\bar{X}^n)\big) > a \Big]
	& \le \mathbb{E} \Big[ \mathbb{P} \big( F(\bar{X}^n) - \langle P^{x}, F \rangle > \frac{a}{2} \, \big \vert \, \bar{X}^n(0) = x \big) \Big]
	\\
	& \qquad + \mathbb{P} \Big[ \langle P^{\bar{X}^n(0)}, F \rangle - \mathbb{E} \big[\langle P^{\bar{X}^n(0)}, F \rangle \big] > \frac{a}{2} \Big]
	\\
	& \le \exp\Big(-\frac{a^2}{8c_1}\Big) + \exp\Big(-\frac{a^2}{8\kappa c_2^2}\Big).
	\end{align*}
	The result \eqref{eq : concentration ineq X bar} follows by choosing $1 / \delta = 8 \max(c_1, \kappa c_2^2)$.
\hfill $\qed$

\medskip

\subsection{Proofs of results in Section~\ref{subsec : concentration toward GS}}	\label{subsec : concentration toward GS proof}

\subsubsection{Proof of Theorem~\ref{thm : W1 L bar}}
	Note that	
	\begin{equation*}
	Y \mapsto \sup_{0 \le t \le T} W_{1} \Big(\frac{1}{n}\sum_{i=1}^n \delta_{Y_i(t)}, \, \widetilde{\mu}_t \Big)
	\end{equation*}
	is $(1/n)$-Lipschitz from $\Big(\big( C ([0, T] : \mathbb{R}^d) \big)^n, \, \Vert\cdot\Vert_{n, 1} \Big)$ to $\mathbb{R}$. Then, for any $a > 0$, 
	\begin{align*}
	\mathbb{P} \big[ \sup_{0 \le t \le T} W_{1} (\bar{L}_{n, t}, \widetilde{\mu}_t) > a \big] 
	&\le \mathbb{P} \Big[ \sup_{0 \le t \le T} W_{1} (\bar{L}_{n, t}, \widetilde{\mu}_t) - \mathbb{E}\big[\sup_{0 \le t \le T} W_{1} (\bar{L}_{n, t}, \widetilde{\mu}_t) \big]  > \frac{a}{2} \Big]
	\\
	& \qquad + \mathbb{P} \Big[ \mathbb{E}\big[\sup_{0 \le t \le T} W_{1} (\bar{L}_{n, t}, \widetilde{\mu}_t) \big]  > \frac{a}{2} \Big].
	\end{align*}
	The first term is bounded by the right-hand side of \eqref{ineq : Wasserstein 2 bound for bar L1} from Theorem~\ref{thm : concen ineq X bar1}.
	
	\smallskip
	
	Let us consider the auxiliary particle system \eqref{eq : SDE} satisfying Assumption~\ref{assump : np(n)}. Corollary~\ref{cor : expectations converging to zero} shows that the last probability vanishes for all but finitely many $n$ and the result follows.
\hfill $\qed$

\smallskip

\subsubsection{Proof of Theorem~\ref{thm : W1 L}}
	We first prove \eqref{ineq : Wasserstein 1 bound for L1}. From the triangle inequality, we obtain
	\begin{align}	\label{ineq : W1 two bounds}
	\mathbb{P} \Big[ \sup_{0 \le t \le T} W_{1} (L_{n, t}, \widetilde{\mu}_t) > a \Big] 
	\le \mathbb{P} \Big[ \sup_{0 \le t \le T} W_{1} (\bar{L}_{n, t}, \widetilde{\mu}_t) > \frac{a}{2} \Big]
	+ \mathbb{P} \Big[ \sup_{0 \le t \le T} W_{1} (L_{n, t}, \bar{L}_{n, t}) > \frac{a}{2} \Big].
	\end{align}
	In what follows, we compute the bound for the last probability on the right-hand side. For fixed $t \in [0, T]$ and $i \in [n]$, we use the notation \eqref{def : P bar} to obtain
	\begin{align}
	X^n_i(t) &- \bar{X}^n_i(t)
	= \int_0^t \sum_{j=1}^n D^{(n)}_{i, j} \phi\big( X^n_i(s), X^n_j(s) \big) \, ds		\label{eq : diff X X bar}
	\\
	&+\int_0^t \sum_{j=1}^n \bar{P}^{(n)}_{i, j} \Big ( \phi\big(X^n_i(s), X^n_j(s) \big) - \phi\big(\bar{X}^n_i(s), \bar{X}^n_j(s) \big) \Big) + \psi\big(X^n_i(s)\big) - \psi\big(\bar{X}^n_i(s) \big) \, ds.		\nonumber
	\end{align}
	We define $\triangle(t) := \frac{1}{n} \sum_{i=1}^n \vert \vert X^n_i - \bar{X}^n_i \vert \vert_{\star, t}$, then deduce from the continuity of $X^n_i(\cdot)-\bar{X}^n_i(\cdot)$ that there exists $t_i \in [0, t]$ for each $i \in [n]$ satisfying
	\begin{align}
	\triangle(t) 
	= \frac{1}{n} \sum_{i=1}^n \vert X^n_i(t_i) &- \bar{X}^n_i(t_i) \vert
	\le \int_0^{t} \frac{1}{n} \sum_{i,j=1}^n \Big\vert D^{(n)}_{i, j} \mathbbm{1}_{[0, t_i]}(s) \phi\big( X^n_i(s), X^n_j(s) \big) \Big\vert \, ds			\label{Dij}
	\\
	&+\frac{1}{n}\int_0^t \sum_{i, j=1}^n \bar{P}^{(n)}_{i, j} \mathbbm{1}_{[0, t_i]}(s) \Big\vert \phi\big(X^n_i(s), X^n_j(s) \big) - \phi\big(\bar{X}^n_i(s), \bar{X}^n_j(s) \big) \Big\vert \, ds 						\label{phi diff2}
	\\
	&+ \frac{1}{n}\int_0^t \sum_{i=1}^n \mathbbm{1}_{[0, t_i]}(s) \Big\vert \psi\big(X^n_i(s)\big) - \psi\big(\bar{X}^n_i(s) \big) \Big\vert \, ds.
	\label{psi diff2}
	\end{align}
	
	\smallskip
	
	Since each component $\phi_k$ of $\phi$ belongs to the $L^1$-Fourier class, there exists a finite complex measure $m_{\phi_k}$ so that we can write for every $k \in [d]$
	\begin{equation}	\label{phi representation}
	\phi_k \big(X^n_i(s), X^n_j(s) \big) = \int_{\mathbb{R}^{2d}} a^k_i(z, s) b^k_j(z, s) m_{\phi_k}(dz), \qquad z = (z_1, z_2),
	\end{equation}
	for some complex functions $a^k_i, b^k_j$ of the form
	\begin{align*}
	a^k_i(z, s) := \exp \big(2\pi \sqrt{-1} \langle X^n_i(s), z_1 \rangle \big),
	\qquad b^k_j(z, s) := \exp \big(2\pi \sqrt{-1} \langle X^n_j(s), z_2 \rangle \big).
	\end{align*}
	
	\smallskip
	
	Using the representation \eqref{phi representation} with an application of H\"older inequality, the integral of \eqref{Dij} is bounded above by
	\begin{align}
	\frac{\max_{1 \le k \le d} \vert\vert m_{\phi_k} \vert\vert_{TM} t}{n} \sum_{i,j=1}^n & \sum_{k=1}^d \int_0^t \int_{\mathbb{R}^{2d}} \Big\vert D^{(n)}_{i, j} \mathbbm{1}_{[0, t_i]}(s) a^k_i(z, s) b^k_j(z, s) m_{\phi_k}(dz) \Big\vert \, ds	\nonumber
	\\
	&= \frac{Kt}{n} \sum_{k=1}^d  \int_0^t \int_{\mathbb{R}^{2d}} \big\vert \langle \mathbf{a}^k(z, s), \, D^{(n)} \mathbf{b}^k(z, s) \rangle \big\vert \, m_{\phi_k}(dz) \, ds,
	\label{abm}
	\end{align}
	where we defined the complex vectors
	\begin{equation}	\label{def : complex vectors}
	\mathbf{a}^k(z, s) := \Big(\mathbbm{1}_{[0, t_i]}(s) a^k_i(z, s) \Big)_{i \in [n]}, \qquad 
	\mathbf{b}^k(z, s) := \big( b^k_j(z, s) \big)_{j \in [n]}, \qquad
	\text{for each } k \in [d].
	\end{equation}
	Since $\ell^{\infty}$-norms of these vectors are bounded by $1$, decomposing them into real and complex parts gives for each $k \in [d]$
	\begin{equation*}
	\big\vert \langle \mathbf{a}^k(z, s), \, D^{(n)} \mathbf{b}^k(z, s) \rangle \big\vert \le 4 \sup \big\{ \langle \mathbf{x}, D^{(n)}\mathbf{y} \rangle : \mathbf{x}, \mathbf{y} \in [-1, 1]^n \big\} = 4\vert\vert D^{(n)} \vert\vert_{\infty \rightarrow 1}.
	\end{equation*}
	Thus, the right-hand side of \eqref{abm} is bounded above by
	\begin{equation*}
	\frac{Kdt^2}{n} \Big( \max_{1 \le k \le d} \vert\vert m_{\phi_k} \vert\vert_{TM} \Big) \vert\vert D^{(n)} \vert\vert_{\infty \rightarrow 1}.
	\end{equation*}
	For the integrals of \eqref{phi diff2} and \eqref{psi diff2}, we use the Lipschitz continuity of $\phi$ and $\psi$, thus there exists a constant $K > 0$ such that
	\begin{equation}		\label{triangle bound}
	\triangle(t) \le K \int_0^t \triangle (s) ds + \frac{\vert\vert D^{(n)} \vert\vert_{\infty \rightarrow 1}}{n}Kt^2, \qquad \forall \, t \in [0, T].
	\end{equation}
	Gr\"onwall's inequality yields
	\begin{equation}	\label{triangle D}
	\triangle(T) \le \frac{K \Vert D^{(n)}\Vert_{\infty \rightarrow 1}}{n},
	\end{equation}
	where $K$ is now a positive constant depending on the time horizon $T$. Recalling the notation $\triangle(t)$, we obtain
	\begin{equation*}
	\sup_{0 \le t \le T} W_{1} (L_{n, t}, \bar{L}_{n, t}) 
	\le \triangle(T)
	\le \frac{K \vert\vert D^{(n)} \vert\vert_{\infty \rightarrow 1}}{n},
	\end{equation*}
	and finally Lemma~\ref{lem : cut norm} gives the bound for the last probability of \eqref{ineq : W1 two bounds}
	\begin{equation}	\label{exponential bound L L bar}
	\mathbb{P} \Big[ \sup_{0 \le t \le T} W_{1} (L_{n, t}, \bar{L}_{n, t}) > \frac{a}{2} \Big] \le \exp \Big( - \frac{a^2 n^2p(n)}{8K^2+2aK/3}\Big).
	\end{equation}
	For the first probability on the right-hand side of \eqref{ineq : W1 two bounds}, Theorem~\ref{thm : W1 L bar} yields for every $n \ge N$
	\begin{equation*}
	\mathbb{P} \Big[ \sup_{0 \le t \le T} W_{1} (\bar{L}_{n, t}, \widetilde{\mu}_t) > \frac{a}{2} \Big] \le 2 \exp \Big( - \frac{\delta a^2n}{16} \Big).
	\end{equation*}
	Thanks to Assumption~\ref{assump : np(n)}, by choosing a larger value for $N \in \mathbb{N}$ than the one in Theorem~\ref{thm : W1 L bar}, we can make
	$\exp \big( - \frac{a^2 n^2p(n)}{8K^2+2aK/3}\big) \le \exp \big( - \frac{\delta a^2n}{16}\big)$ for every $n \ge N$, and the assertion \eqref{ineq : Wasserstein 1 bound for L1} follows.
	
	\smallskip
	
	For the result \eqref{ineq : Wasserstein 1 bound for BL}, we can approximate general $\phi$ with those in $L^1$-Fourier class by the approximation method in Section~5.1.3 of \cite{Oliveria:Reis}, to find the exponential bound for the probability $\mathbb{P} [ \sup_{0 \le t \le T} d_{BL} (L_{n, t}, \bar{L}_{n, t}) > a/2 ]$ similar to \eqref{exponential bound L L bar}. By recalling the fact $d_{BL} \le W_1$ and replacing all the $W_1$-metrics with the $d_{BL}$-metrics in \eqref{ineq : W1 two bounds}, we arrive at \eqref{ineq : Wasserstein 1 bound for BL}.
\hfill $\qed$

\smallskip

\subsubsection{Proof of Theorem~\ref{thm : W2 L bar}}
	It is easy to verify the $(\frac{1}{\sqrt{n}})$-Lipschitz continuity of the map 	
	\begin{equation*}
	Y \mapsto \sup_{0 \le t \le T} W_{2} \Big(\frac{1}{n}\sum_{i=1}^n \delta_{Y_i(t)}, \, \widetilde{\mu}_t \Big),
	\end{equation*}
	from $\Big(\big( C ([0, T] : \mathbb{R}^d) \big)^n, \, \Vert\cdot\Vert_{n, 2} \Big)$ to $\mathbb{R}$. Then, for any $a > 0$, 
	\begin{align*}
	\mathbb{P} \big[ \sup_{0 \le t \le T} W_2 (\bar{L}_{n, t}, \widetilde{\mu}_t) > a \big]
	&\le \mathbb{P} \Big[ \sup_{0 \le t \le T} W_2 (\bar{L}_{n, t}, \widetilde{\mu}_t) - \mathbb{E}\big[\sup_{0 \le t \le T} W_2 (\bar{L}_{n, t}, \widetilde{\mu}_t) \big]  > \frac{a}{2} \Big]
	\\
	& \qquad + \mathbb{P} \Big[ \mathbb{E}\big[\sup_{0 \le t \le T} W_2 (\bar{L}_{n, t}, \widetilde{\mu}_t) \big]  > \frac{a}{2} \Big].
	\end{align*}
	The first term is bounded by the right-hand side of \eqref{ineq : Wasserstein 2 bound for bar L} from Theorem~\ref{thm : concen ineq X bar2}. As in the proof of Theorem~\ref{thm : W1 L bar}, the last probability vanishes for all but finitely many $n$ from Corollary~\ref{cor : expectations converging to zero}.
\hfill $\qed$

\smallskip

\subsubsection{Proof of Lemma~\ref{lem : DTD}}
	We note from \eqref{def : P bar} that $\{D^{(n)}_{i, j}\}_{1 \le i, j \le n}$ are independent zero-mean random variables and for every $i, j \in [n]$
	\begin{equation}	\label{eq : expectation of D^2}
		\mathbb{E}\big[(D^{(n)}_{i, j})^2\big] = \frac{p(n)G(\frac{i}{n}, \frac{j}{n}) \big(1 - p(n)G(\frac{i}{n}, \frac{j}{n}) \big)}{\big(np(n)\big)^2}.
	\end{equation}
	In particular, since $p(n) \le 1$, we have $0 \le p(n)G(\frac{i}{n}, \frac{j}{n}) \le 1$, and thus $\mathbb{E}\big[(D^{(n)}_{i, j})^2\big] \le 1/(4n^2p(n)^2)$.
	
	\smallskip
	
	Let us fix any $n \in \mathbb{N}$. For arbitrary $n$-dimensional vectors $\mathbf{x}, \, \mathbf{y} \in [-1, 1]^n$, we have
	\begin{align*}
	& \quad \langle \mathbf{x}, (D^{(n)})^\top D^{(n)}\mathbf{y} \rangle 
	= \sum_{i=1}^n \big(D^{(n)} \mathbf{x}\big)_i \big(D^{(n)} \mathbf{y} \big)_i
	= \sum_{i=1}^n \sum_{j=1}^n \sum_{\substack{k=1 \\ k \neq j}}^n D^{(n)}_{i, j} D^{(n)}_{i, k} x_j y_k 
	+ \sum_{i=1}^n \sum_{j=1}^n \big( D^{(n)}_{i, j} \big)^2 x_j y_j
	\\
	&\le \sum_{i=1}^n \sum_{j=1}^n \sum_{\substack{k=1 \\ k \neq j}}^n D^{(n)}_{i, j} D^{(n)}_{i, k} x_j y_k
	+ \sum_{i=1}^n \sum_{j=1}^n \Big( \big( D^{(n)}_{i, j} \big)^2 x_j y_j - \mathbb{E}\big[ \big(D^{(n)}_{i, j}\big)^2 x_j y_j \big] \Big)
	+ \sum_{i=1}^n \sum_{j=1}^n \frac{x_i y_j}{4n^2p(n)^2}.
	\end{align*}
	Thus, we have for fixed arbitrary $\eta > 0$
	\begin{align}
		\mathbb{P} \bigg[ \frac{\Vert (D^{(n)})^\top D^{(n)} \Vert_{\infty \rightarrow 1}}{n} > \eta \bigg] 
		&\le \mathbb{P} \bigg[ \frac{1}{n} \sum_{i=1}^n \sum_{j=1}^n \sum_{\substack{k=1 \\ k \neq j}}^n D^{(n)}_{i, j} D^{(n)}_{i, k} x_j y_k > \frac{\eta}{3} \bigg]		\nonumber
		\\
		&+ \mathbb{P} \bigg[ \frac{1}{n}\sum_{i=1}^n \sum_{j=1}^n \Big( \big( D^{(n)}_{i, j} \big)^2 x_j y_j - \mathbb{E}\big[ \big(D^{(n)}_{i, j}\big)^2 x_j y_j \big] \Big) > \frac{\eta}{3} \bigg]		\nonumber
		\\
		&+ \mathbb{P} \bigg[ \frac{1}{4np(n)^2} > \frac{\eta}{3} \bigg]
		=: P_1 + P_2 + P_3.									\label{P1P2P3}
	\end{align}
	From Assumption~\ref{assump : np(n)^2}, there exists $N \in \mathbb{N}$ such that $P_3$ vanishes for every $n \ge N$. In the following, we find the bounds for $P_1$ and $P_2$. Using the distribution
	\begin{equation*}
		D^{(n)}_{i, j} =
		\begin{cases}
			\frac{1-p(n) G(i/n, j/n)}{np(n)} \le \frac{1}{np(n)},	&\text{ with probability } p(n)G(i/n, j/n),
			\\
			\\
			-\frac{1}{n}G(i/n, j/n) \ge - \frac{1}{n},	&\text{ with probability } 1-p(n)G(i/n, j/n),
		\end{cases}
	\end{equation*}
	for each $i, j \in [n]$, we derive for $P_1$
	\begin{align*}
		P_1 &\le \sum_{i=1}^n \mathbb{P} \bigg[ \sum_{j=1}^n \sum_{\substack{k=1 \\ k \neq j}}^n D^{(n)}_{i, j} D^{(n)}_{i, k} x_j y_k > \frac{\eta}{3} \bigg]
		\le \sum_{i=1}^n \sum_{j=1}^n \mathbb{P} \bigg[ \sum_{\substack{k=1 \\ k \neq j}}^n D^{(n)}_{i, j} D^{(n)}_{i, k} x_j y_k > \frac{\eta}{3n} \bigg]
		\\
		& = \sum_{i=1}^n \sum_{j=1}^n \mathbb{E} \bigg[ \mathbb{P} \Big( \sum_{\substack{k=1 \\ k \neq j}}^n D^{(n)}_{i, j} D^{(n)}_{i, k} x_j y_k > \frac{\eta}{3n} \Big\vert D^{(n)}_{i, j} \Big) \bigg]
		\\
		& \le \sum_{i=1}^n \sum_{j=1}^n \bigg( \mathbb{P} \Big[ \sum_{\substack{k=1 \\ k \neq j}}^n D^{(n)}_{i, k} x_j y_k > \frac{\eta p(n)}{3} \Big] + \mathbb{P} \Big[ \sum_{\substack{k=1 \\ k \neq j}}^n D^{(n)}_{i, k} x_j y_k < -\frac{\eta}{3} \Big] \bigg).
	\end{align*}
	The summands $D^{(n)}_{i, k} x_j y_k$ in the last two probabilities are independent zero-mean random variables, bounded above by $1/(np(n))$, bounded below by $-1/n$, and satisfies
	\begin{equation*}
		\sum_{\substack{k=1 \\ k \neq j}}^n \mathbb{E} \big[ (D^{(n)}_{i, k} x_j y_k)^2 \big] \le \frac{n-1}{4n^2p(n)^2} \le \frac{1}{4np(n)^2},
	\end{equation*}
	from \eqref{eq : expectation of D^2}. From Bernstein's inequality~(Lemma~\ref{lem : Bernstein}), we have
	\begin{align*}
		\mathbb{P} \Big[ \sum_{\substack{k=1 \\ k \neq j}}^n D^{(n)}_{i, k} x_j y_k > \frac{\eta p(n)}{3} \Big] 
		&\le \exp \bigg( - \frac{2np(n)^4 \eta^2}{9+ 4\eta p(n)^2} \bigg),
		\\
		\mathbb{P} \Big[ -\sum_{\substack{k=1 \\ k \neq j}}^n D^{(n)}_{i, k} x_j y_k > \frac{\eta}{3} \Big] 
		&\le \exp \bigg( - \frac{2np(n)^2 \eta^2}{9+ 4\eta p(n)^2} \bigg),
	\end{align*}
	thus
	\begin{equation}	\label{P1}
		P_1 \le 2n^2 \exp \bigg( - \frac{2np(n)^4 \eta^2}{9+ 4\eta p(n)^2} \bigg).
	\end{equation}
	
	\smallskip
	
	We now compute the bound for $P_2$. Since we have
	\begin{equation}	\label{P2 bound}
		P_2 \le \sum_{i=1}^n \mathbb{P} \bigg[ \sum_{j=1}^n \Big( \big( D^{(n)}_{i, j} \big)^2 x_j y_j - \mathbb{E}\big[ \big(D^{(n)}_{i, j}\big)^2 x_j y_j \big] \Big) > \frac{\eta}{3} \bigg]
	\end{equation}
	and the summands $\big( D^{(n)}_{i, j} \big)^2 x_j y_j - \mathbb{E}\big[ \big(D^{(n)}_{i, j}\big)^2 x_j y_j \big]$ in the probability are independent zero-mean random variables bounded above by $5/(2np(n))^2$. Moreover, we easily obtain the bound $\mathbb{E}\big[ (D^{(n)}_{i, j})^4\big] \le 1/(np(n))^4$, and thus the sum of variances of summands are
	\begin{align*}
		\sum_{j=1}^n \mathbb{E} \bigg[ \Big( \big( D^{(n)}_{i, j} \big)^2 x_j y_j - \mathbb{E}\big[ \big(D^{(n)}_{i, j}\big)^2 x_j y_j \big] \Big)^2 \bigg]
		&\le \sum_{j=1}^n \mathbb{E} \bigg[ \Big( \big( D^{(n)}_{i, j} \big)^2 - \mathbb{E}\big[ \big(D^{(n)}_{i, j}\big)^2 \big] \Big)^2 \bigg]
		\\
		& \le \sum_{j=1}^n \bigg( \mathbb{E} \big[ (D^{(n)}_{i, j})^4\big] + \Big(\mathbb{E} \big[ (D^{(n)}_{i, j})^2\big] \Big)^2 \bigg)
		\le \frac{17}{16n^3p(n)^4}.
	\end{align*}
	Applying Bernstein's inequality~(Lemma~\ref{lem : Bernstein}) to each probability of \eqref{P2 bound} yields
	\begin{equation}	\label{P2}
		P_2 \le n \exp \bigg( - \frac{8n^3p(n)^4 \eta^2}{153+ 20 np(n)^2\eta} \bigg).
	\end{equation}
	Comparing the bounds of \eqref{P1}, \eqref{P2}, modifying the value of $N$ if necessary, and plugging into \eqref{P1P2P3}, the result follows.	
\hfill $\qed$

\smallskip

\subsubsection{Proof of Theorem~\ref{thm : W2 L}}
	The argument is similar to the proof of Theorem~\ref{thm : W1 L}. The triangle inequality gives
	\begin{align}	\label{ineq : W2 two bounds}
	\mathbb{P} \Big[ \sup_{0 \le t \le T} W_{2} (L_{n, t}, \widetilde{\mu}_t) > a \Big] 
	\le \mathbb{P} \Big[ \sup_{0 \le t \le T} W_{2} (\bar{L}_{n, t}, \widetilde{\mu}_t) > \frac{a}{2} \Big]
	+ \mathbb{P} \Big[ \sup_{0 \le t \le T} W_{2} (L_{n, t}, \bar{L}_{n, t}) > \frac{a}{2} \Big].
	\end{align}
	Recalling the identity \eqref{eq : diff X X bar}, applying a series of H\"older's inequality, and using the Lipschitz property yield
	\begin{align*}
		\big\vert X^n_i(t) &- \bar{X}^n_i(t) \big\vert^2
		- 2t \int_0^t \big\vert \sum_{j=1}^n D^{(n)}_{i, j} \phi\big( X^n_i(s), X^n_j(s) \big) \Big\vert^2 \, ds
		\\
		& \le 2K^2t \int_0^t 2 \bigg( \sum_{j=1}^n \bar{P}^{(n)}_{i, j} \Big( \big\vert X^n_i(s)-\bar{X}^n_i(s) \big\vert + \big\vert X^n_j(s)-\bar{X}^n_j(s) \big\vert \Big)^2 + 2\big\vert X^n_i(s)-\bar{X}^n_i(s) \big\vert^2 \, ds
		\\
		& \le 2K^2t \int_0^t 6 \big\vert X^n_i(s)-\bar{X}^n_i(s) \big\vert^2 + \frac{2}{n} \sum_{j=1}^n \big\vert X^n_j(s)-\bar{X}^n_j(s) \big\vert^2 \, ds, \qquad \forall \, t \in [0, T].
	\end{align*}
	For a fixed $t \in [0, T]$, by using the continuity of $X^n_i(\cdot)-\bar{X}^n_i(\cdot)$, there exists $t_i \in [0, t]$ for each $i \in [n]$ satisfying $\Box(t) := \frac{1}{n} \sum_{i=1}^n \Vert X^n_i - \bar{X}^n_i \Vert^2_{\star, t} = \frac{1}{n} \sum_{i=1}^n \big\vert X^n_i(t_i) - \bar{X}^n_i(t_i) \big\vert^2$. Combining with the last inequality, we have
	\begin{equation}	\label{ineq : Box}
	\Box(t)
	\le 2t \int_0^{t} \frac{1}{n} \sum_{i=1}^n \Big\vert \sum_{j=1}^n D^{(n)}_{i, j} \mathbbm{1}_{[0, t_i]}(s) \phi\big( X^n_i(s), X^n_j(s) \big) \Big\vert^2 \, ds
	+16K^2t \int_0^t \Box(s) \, ds.
	\end{equation}
	We recall the representations \eqref{phi representation} and \eqref{def : complex vectors} and use H\"older's inequality to derive for the first integral on the right-hand side
	\begin{align*}
		\int_0^{t} \frac{1}{n} \sum_{i=1}^n \Big\vert \sum_{j=1}^n D^{(n)}_{i, j} &\mathbbm{1}_{[0, t_i]}(s) \phi\big( X^n_i(s), X^n_j(s) \big) \Big\vert^2 \, ds
		\\
		&\le d \int_0^{t} \frac{1}{n} \sum_{i=1}^n \sum_{k=1}^d \bigg\vert \sum_{j=1}^n \int_{\mathbb{R}^{2d}} \big\vert D^{(n)}_{i, j} \mathbbm{1}_{[0, t_i]}(s) a^k_i(z, s) b^k_j(z, s) \big\vert \, m_{\phi_k}(dz) \bigg\vert^2 \, ds
		\\
		& \le d \int_0^{t} \frac{1}{n} \sum_{i=1}^n \sum_{k=1}^d \Vert m_{\phi_k} \Vert_{TM} \int_{\mathbb{R}^{2d}} \big\vert \big(\mathbf{a}^k(z, s)\big)_i \, \big( D^{(n)} \mathbf{b}^k(z, s)\big)_i \big\vert^2 \, m_{\phi_k}(dz) \, ds
		\\
		& \le d \int_0^{t} \sum_{k=1}^d \Vert m_{\phi_k} \Vert_{TM} \int_{\mathbb{R}^{2d}} \frac{1}{n} \sum_{i=1}^n \big\vert \big( D^{(n)} \mathbf{b}^k(z, s)\big)_i \big\vert^2 \, m_{\phi_k}(dz) \, ds
		\\
		& \le \max_{1 \le k \le d} \Vert m_{\phi_k} \Vert_{TM}^2 d^2 t \frac{\Vert (D^{(n)})^\top D^{(n)} \Vert_{\infty \rightarrow 1}}{n}.
	\end{align*}
	In the last two inequalities, we used the fact that $\ell^{\infty}$-norms of two vectors $\mathbf{a}^k(z, s)$ and $\mathbf{b}^k(z, s)$ are bounded by $1$. Thus, from \eqref{ineq : Box}, there exists a constant $K>0$ such that
	\begin{equation*}
		\Box(t) \le Kt^2 \frac{\Vert (D^{(n)})^\top D^{(n)} \Vert_{\infty \rightarrow 1}}{n} + Kt \int_0^t \Box(s) ds
	\end{equation*}
	holds for every $t \in [0, T]$ and applying Gr\"onwall's inequality gives
	\begin{equation*}
		\Box(T) \le \frac{K\Vert (D^{(n)})^\top D^{(n)} \Vert_{\infty \rightarrow 1}}{n},
	\end{equation*}
	where the constant $K$ now depends on $T$.
	
	\smallskip
	
	Since we have
	\begin{equation*}
		\sup_{0 \le t \le T} W_2(L_{n, t}, \bar{L}_{n, t}) \le \sup_{0 \le t \le T} \sqrt{\frac{1}{n} \sum_{i=1}^n \big\vert X^n_i(t) - \bar{X}^n_i(t) \big\vert^2}
		\le \sqrt{\Box(T)}, 
	\end{equation*}
	Lemma~\ref{lem : DTD} shows that there exists $N \in \mathbb{N}$ such that the last probability in \eqref{ineq : W2 two bounds} has the bound
	\begin{align*}
		\mathbb{P} \Big[ \sup_{0 \le t \le T} W_{2} (L_{n, t}, \bar{L}_{n, t}) > \frac{a}{2} \Big]
		&\le \mathbb{P} \Big[ \Box(T) > \frac{a^2}{4} \Big]
		\le \mathbb{P} \Big[ \frac{\Vert (D^{(n)})^\top D^{(n)} \Vert_{\infty \rightarrow 1}}{n} > \frac{a^2}{4K} \Big]
		\\
		&\le 3n^2 \exp \bigg( - \frac{a^4 np(n)^4}{72K^2+8a^2K} \bigg),
	\end{align*}
	for every $n \ge N$.
	
	\smallskip
	
	On the other hand, Theorem~\ref{thm : W2 L bar} gives the bound for the other probability in \eqref{ineq : W2 two bounds}. By comparing both of the bounds under Assumption~\ref{assump : np(n)^2}(a), the assertion \eqref{ineq : Wasserstein 2 bound for L1} follows. The result \eqref{ineq : Wasserstein 2 bound for L1 exponential} is now clear under Assumption~\ref{assump : np(n)^2}(b), by setting $p(n) \equiv 1$ and redefining the constants $\delta > 0$ and $N \in \mathbb{N}$ appropriately.
\hfill $\qed$

\bigskip

\bigskip

\bigskip


\renewcommand{\bibname}{References}
\bibliography{aa_bib}
\bibliographystyle{apalike}

\end{document}